\documentclass[reqno,11pt]{amsart}

\usepackage{color}

\newtheorem{theorem}{Theorem}[section] 
\newtheorem{lemma}[theorem]{Lemma}
\newtheorem{corollary}[theorem]{Corollary}
 \theoremstyle{remark}
\newtheorem{remark}[theorem]{Remark}

\newcommand\cD{\mathcal{D}}
\newcommand\cL{\mathcal{L}}

\newcommand\bL{\mathbb{L}}
\newcommand\bH{\mathbb{H}}
\newcommand\bR{\mathbb{R}}

\newcommand\Lie{{\rm Lie\,}}

\renewcommand\det{{\rm det}\,}

\newcommand{\mysection}[1]{\section{#1}
      \setcounter{equation}{0}}

\newcommand{\nlimsup}{\operatornamewithlimits{\overline{lim}}}

\theoremstyle{definition}
\newtheorem{assumption}{Assumption}[section]
\newtheorem{definition}{Definition}[section]

\newcommand\cbrk{\text{$]$\kern-.15em$]$}}
\newcommand\opar{\text{\raise.2ex\hbox{${\scriptstyle | }$}\kern-.34em$($} }

 \makeatletter
 \def\dashint{%
 \operatorname%
 {\,\,\text{\bf--}\kern-.98em\DOTSI\intop\ilimits@\!\!}}
 \makeatother

 \newcommand{\WO}{\overset{\scriptscriptstyle0}%
{ W}\,\!}

\newcommand{\cHO}{\overset{\,\,\scriptscriptstyle0}%
{\mathcal H}\,\!}

\renewcommand\div{\text{\rm div}\,}

\newcommand\cF{\mathcal{F}}

\newcommand\cZ{\mathcal{Z}}
\newcommand\cH{\mathcal{H}}

\begin{document}

\title[Hypoellipticity
for filtering problems]{Hypoellipticity
for filtering problems of partially observable diffusion processes} 
\author{N.V. Krylov}
\thanks{The  author was partially supported by
 NSF Grant DMS-1160569}
\email{krylov@math.umn.edu}
\address{127 Vincent Hall, University of Minnesota,
 Minneapolis, MN, 55455}

\begin{abstract}
We prove that under H\"ormander's type conditions
on the coefficients of the unobservable component
of a partially observable
diffusion process
the filtering density  is infinitely differentiable
and can be represented as the integral of
an infinitely differentiable kernel against the
  prior initial distribution.
These results are derived from more general results
obtained for SPDEs. One the main novelty
of the paper is the existence and smoothness of the kernel,
another that we allow the coefficients of our
partially observable process to be just measurable with
respect to the time variable.
 
\end{abstract}



\maketitle

\mysection{Introduction}
                                          \label{section 6.28.1}

Let $(\Omega,\cF,P)$ be a complete probability space
with an increasing filtration $\{\cF_{t},t\geq0\}$
of complete with respect to $(\cF,P)$ $\sigma$-fields
$\cF_{t}\subset\cF$.   Let
 $w^{k}_{t}$, $k=1,2,...,d_{1}$, be independent one-dimensional
Wiener processes with respect to $\{\cF_{t}\}$, where $d_{1}\geq1$
is an integer.  

Let   $d,d'\geq1$  
be integers.
Consider a $d+d'$-dimensional two-component process
 $z_{t}=(x_{t},y_{t})$
with $x_{t}$ being $d$-dimensional and $y_{t}$ 
$d' $-dimensional.   We
assume that $z_{t}$ is a diffusion process
 defined as a solution of the system
\begin{equation}\begin{split}         
                                            \label{eq3.2.14} 
& dx_{t}=b(t,z_{t}) dt+\theta^{k} (t,z_{t})\,dw^{k}_{t}, \\ 
& dy_{t}=B(t,z_{t}) dt+\Theta^{k}(t,y_{t})\,dw^{k}_{t}
\end{split}
\end{equation}
with some initial data independent of the process $w_{t}$.
The coefficients of \eqref{eq3.2.14} are assumed to be
 vector-valued
functions of appropriate dimensions defined on
 $[0,\infty)\times\bR^{d+d'}$.
Actually $\Theta^{k}(t,y)$ are assumed to be independent
 of $x$, so  that they are
  functions on $[0,\infty)\times\bR^{d' }$ rather than 
$[0,\infty)\times\bR^{d+d'}$
but as always we may think of $\Theta^{k}(t,y)$ as  
 functions of $(t,z)$ as well.

One of the main goals of the paper is to show that under
H\"ormander's type conditions satisfied for $x$
lying in a ball $B$, in some sense uniformly with respect to
$t$ and $y$, there exists a  
function $p(t,y,x)=p(\omega,t,y,x)\geq0$, which is infinitely
differentiable with respect to $(y,x)\in B^{2}$
for any $t>0$ and $\omega$, such that for any 
$f\in C^{\infty}_{0}(B)$ and $t>0$ with probability one
$$
E\{f(x_{t})\mid \cF^{y}_{t}\}=
\int_{B}\int_{B}f(x)p(t,y,x)\,P_{0}(dy)dx,
$$
where $P_{0}$ is the conditional distribution
of $x_{0}$ given $y_{0}$ and $\cF^{y}_{t}
=\sigma\{y_{s},s\leq t\}$. Naturally,
\begin{equation}
                                            \label{9.10.1}
\int_{B} p(t,y,x)\,P_{0}(dy)
\end{equation}
turns out to be infinitely differentiable with respect
to $x\in B$ and represent the conditional density  $\pi_{t}(x) $
of $x_{t}$ given $\cF^{y}_{t}$. 

In the literature two approaches to prove
infinite differentiability of $\pi_{t}(x)$
for degenerate processes
under H\"ormander's type conditions are known.
The first one is based on filtering equations for $\pi_{t}$,
which are stochastic partial differential equations (SPDEs).
This approach was initiated by Wentzell \cite{We65} and
in  a more general and {\em time inhomogeneous} case
outlined by Kunita in \cite{Ku81} and \cite{Ku82}.
It is worth noting that in \cite{We65} the coefficient 
$B$ is supposed to be independent of $x$ and in \cite{Ku81}
the functions $b$ and $\theta$ are independent of $y$.
Equations in \cite{Ku82} seem not to cover general 
filtering equations either.
In \cite{We65}, \cite{Ku81}, and \cite{Ku82} the SPDE
  is reduced to an ordinary  parabolic
equation with random coefficients by using
a random change of coordinates.
Without this reduction 
 Chaleyat-Maurel and Michel in \cite{CM} achieve the goal
in the {\em time homogeneous case}   by
mimicking some steps which are used in the proof of 
the deterministic
H\"ormander theorem. In their case as well as in \cite{Ch}
the matrix $(\Theta^{1},...,\Theta^{d'})$ is assumed to have form
$(0,I)$ where $I$ is $d'\times d'$ identity matrix.

However, some of the arguments in \cite{Ku81} and \cite{Ku82} are
based on the claim that H\"ormander's type theorem holds
and can be proved
by using Malliavin calculus for equations whose coefficients
are only continuous with respect to $t$. Such a proof is unknown
even now. It also looks like in \cite{CM} there is a gap
at the point  where the authors claim without proof that,
roughly speaking,
what holds for the unknown function also holds
for its fractional derivatives. 

Our approach is also based on using filtering equations
but since we allow the coefficients of \eqref{eq3.2.14} 
to be just measurable with respect to $t$ our type
of H\"ormander's condition is more restrictive than
in \cite{CM} where the coefficients are time independent.
In contrast with \cite{Ku81} and \cite{Ku82} we do not 
appeal to 
Malliavin calculus and instead rely on some analytical facts
which we prove for more general SPDEs.

The second approach to proving
infinite differentiability of $\pi_{t}(x)$ almost
completely ignores filtering equations and
is based on Malliavin calculus and first appeared
in the paper by Bismut and  Michel \cite{BM}.
In their model the coefficients are time independent,
but the H\"ormander type condition imposed, albeit global,
is much weaker than ours.
 Kusuoka and Stroock \cite{KS} further relax
the H\"ormander type condition in \cite{BM}
again in time independent case
but in what concerns filtering they assume that 
$B(t,x,y)$ is independent of $x$, so that the problem
 becomes a problem in the theory
  of conditional Markov processes
rather then a more or less general filtering problem,
 because the coefficients of the
equation for the observation process $y_{t}$ are supposed to be {\em independent}
of the signal process $x_{t}$.
This result
can also be found in \cite{Ro}.
In the recent publication by Chaleyat-Maurel \cite{Ch}
and references therein
one can find a detailed account of the progress
concerning Malliavin calculus and filtering equations.
In particular, in \cite{Ch} local versions of
  H\"ormander's type condition from \cite{CM} is used to obtain
the local regularity of solutions.
It seems that these methods are not applicable in
our case of coefficients only
 measurable with respect to $t$.

Apart from this novelty concerning time dependence,
the fact that
the conditional density $\pi_{t}(x) $
 is represented as  \eqref{9.10.1}
with infinitely differentiable kernel seems to be new
for degenerate diffusions $z_{t}$ under H\"ormander's
type condition.

We derive our results about filtering densities
in Section \ref{section 9.10.1} from results
of Sections \ref{section 7.16.1}, \ref{section 8.26.1},
and \ref{section 9.10.2}. In these sections we treat
more general SPDEs than the filtering equation.

The reader understands that  $\bR^{d}$
is a Euclidean space of column-vectors (written in a common
abuse of notation as) $x=(x^{1},...,x^{d})$.
Denote 
$$
D_{i}=\partial/\partial x^{i},\quad D_{ij}=D_{i}D_{j},
$$
and 
for an $\bR^{d}$-valued function 
$\sigma_{t}(x)=\sigma_{t}(\omega,x)$ on $\Omega\times[0,
\infty)\times\bR^{d}$ 
and functions $u_{t}(x)=u_{t}(\omega,x)$
on $\Omega\times[0,\infty)\times\bR^{d}$ set
$$
L_{\sigma_{t}}u_{t}(x)=[D_{i}u_{t}(x)]\sigma^{i}_{t}(x).
$$

Next take an integer $d_{2}\geq 1$ and 
assume that we are given $\bR^{d}$-valued functions
$\sigma^{k}_{t}=(\sigma^{ik}_{t})$, $k=0,...,d_{2}+d_{1}$, on 
$\Omega\times[0,\infty)\times\bR^{d}$,
which are infinitely differentiable with respect to $x$
for any $(\omega,t)$, and define the 
operator 
\begin{equation}
                                                       \label{7.22.4}
L_{t}=(1/2)\sum_{k=1}^{d_{2}+d_{1}}L^{2}_{\sigma^{k}_{t}}
+L_{\sigma^{0}_{t}}.
\end{equation}
Assume that on $\Omega\times[0,\infty)\times\bR^{d}$
we are also given certain real-valued infinitely differentiable
functions $c_{t}(x)$ and $\nu^{k}_{t}(x)$, $k=1,...,d_{1}$,
and that on $\Omega\times[0,\infty)\times
\bR^{d}$ we are given  real-valued
functions $f_{t}$ and $g_{t}^{k}$,  $k=1,...,d_{1}$. Then under natural
additional assumptions which will be specified later the
SPDE
\begin{equation}
                                               \label{6.28.1}
du_{t}=(L_{t}u_{t}+c_{t}u_{t}+f_{t})\,dt
+(L_{\sigma^{ k}_{t}}u_{t}+\nu^{k}_{t}u_{t}+g^{k}_{t})\,dw^{k}_{t}
\end{equation}
makes sense (where and below the summation convention over
repeated indices is enforced regardless of whether they
stand at the same level or at different ones).

One of the main results
 of this paper is Theorem \ref{theorem 7.24.1} saying that
if the initial condition is a generalized functions
of class $H^{-n}_{2}$, then \eqref{6.28.1} has a unique
solution with this initial data without any nondegeneracy
or H\"ormander's type condition. Before this result was known
only if $n\geq1$ is an integer (see \cite{KR_82}). The result is
important because it allows one to take a $\delta$-function
as the initial condition.

After the existence of solutions is secured we continue our
investigation under H\"ormander's type condition and
in Section \ref{section 8.26.1}
prove, roughly speaking,
that, if   $ (s_{1},s_{2})\in(0,T)$ and for any $ \omega\in\Omega$
and $t\in(s_{1},s_{2})$
the Lie algebra generated by the vector-fields $\sigma^{d_{1}+k}_{t}$,
$k=1,...,d_{2}$, has dimension $d$ everywhere in $B_{R}$
and $f_{t}$ and $g_{t}^{k}$ are infinitely differentiable
in $B_{R}$ for any $ \omega\in\Omega$
and $t\in(s_{1},s_{2})$, then the generalized function
$u_{t}$ satisfying \eqref{6.28.1} coincides on 
$(s_{1},s_{2})\times B_{R}$ with a function which is infinitely
differentiable with respect to $x$.
In Section \ref{section 9.10.1} we apply this result
to filtering problems. In the same section
we apply the results of Section \ref{section 9.10.2}
to derive the existence of smooth filtering kernels.
The results of Section \ref{section 9.10.2}
bear on kernels (or fundamental solutions)
for more general SPDEs.

 In the whole article $T$ is a fixed number from $(0,\infty)$.

\mysection{An existence theorem for SPDEs}
                                    \label{section 7.16.1} 

Denote by $\cD$  the space 
of
generalized functions on  $\bR^{d}$, and as usual introduce
$\Lambda=(1-\Delta)^{1/2}$, $H^{n}_{2}=\Lambda^{-n}\cL_{2}$,
where $\cL_{2}$ is the Hilbert
 space of real-valued square integrable
functions on $\bR^{d}$ with usual norm. The scalar product
in $H^{n}_{2}$ will be denoted by $(\cdot,\cdot)_{n}$.

Denote by $H^{n}_{2}(T)$ the set of $\cD$-valued functions
$u_{t}$ on $\Omega\times[0,T]$ such that $(u_{t},\phi)$
is predictable and $(u_{0},\phi)$ is $\cF_{0}$-measurable
for any $\phi\in C^{\infty}_{0}(\bR^{d})$
 and
$$
\int_{0}^{T}\|u_{t}\|^{2}_{n}\,dt<\infty
\quad \text{(a.s.)}.
$$
Sometimes it is necessary to indicate   which
filtration of $\sigma$-fields is involved in the
definition of predictable functions. In these
cases we write $H^{n}_{2}(T)=H^{n}_{2}(T, \cF_{\cdot})$

Introduce $\cH^{n}_{2}$ as the set of
$\cF_{0}$-measurable $H^{n}_{2}$-valued
function on $\Omega$. For an open ball $B$
by $\cHO^{n}_{2}(B)$ we mean the subset of $\cH^{n}_{2}$
consisting of generalized functions with 
(closed) support in $B$.
Define
$$
 \bH^{n}_{2}(T)=\{u\in H^{n}_{2}(T):
E\int_{0}^{T}\|u_{t}\|^{2}_{n}\,dt<\infty\}.
$$

\begin{assumption}
                                      \label{assumption 6.28.1}
(i)  The  functions $\sigma^{k}_{t}(x)$, $k=0,...,d_{2}+d_{1}$,
$c_{t}$, $\nu^{k}_{t}$, $k=1,...,d_{1}$, are infinitely
differentiable with respect to $x$ and each of their derivatives
of any order
is bounded on $\Omega\times[0,T]\times \bR^{d}$.
These functions are predictable with respect to $(\omega,t)$
for any $x\in \bR^{d}$;

(ii) For an $n\in\bR$, we have that 
$ f\in  H^{n}_{2}(T)$, $g^{k}\in  H^{n+1}_{2}(T)$, 
$k=1,...,d_{1}$,
and

(iii)     $u_{0}\in \cH^{n}_{2}$.
\end{assumption}

\begin{definition}
                               \label{definition 7.31.1}
By a {\em normal\/} solution of \eqref{6.28.1}
of class $H^{n}_{2}(T)$  with initial
condition $u_{0}$ we mean a function $u$ which belongs
to $H^{n}_{2}(T)$, such that (a.s.) $u_{t}$
is a continuous $H^{n-1}_{2}$-valued function and
with probability one
\begin{equation}
                                       \label{7.31.1}
u_{t}=u_{0}+\int_{0}^{t}
(L_{s}u_{t}+c_{s}u_{s}+f_{s})\,ds
+\int_{0}^{t}
(L_{\sigma^{ k}_{s}}u_{s}+\nu^{k}_{s}u_{s}+
g^{k}_{s})\,dw^{k}_{s}
\end{equation}
for all $t\in[0,T]$.
\end{definition}

\begin{remark}
                                \label{remark 7.31.1}
The usual and stochastic integrals of
Hilbert space valued functions are well defined,
so that the right-hand side of \eqref{7.31.1}
is a continuous $H^{n-2}_{2}$-valued process.
\end{remark}
\begin{remark}
                                \label{remark 7.31.2}
We  say that a function $u$ of class
$H^{n}_{2}(T)$ is a {\em generalized\/} solution
of \eqref{6.28.1} with initial
condition $u_{0}$ if for any 
$\phi\in C^{\infty}_{0}(\bR^{d})$
$$
(u_{t},\phi)=(u_{0},\phi)+\int_{0}^{t}
(L_{s}u_{t}+c_{s}u_{s}+f_{s},\phi)\,ds
$$
\begin{equation}
                                       \label{7.31.2}
+\int_{0}^{t}(L_{\sigma^{ k}_{s}}u_{s}+\nu^{k}_{s}u_{s}+
g^{k}_{s},\phi)\,dw^{k}_{s}
\end{equation}
for almost all $(\omega,t)\in\Omega\times[0,T]$,
where by $(\cdot,\cdot)$ we mean the pairing between
test functions and generalized ones.

By the way, recall that if $u\in H^{n}_{2}$
and $\phi\in C^{\infty}_{0}(\bR^{d})$, then
$$
(u,\phi)=(\Lambda^{k}u,\Lambda^{m} \phi)_{0},
$$
as long as $k\leq n$ and $k+m=2n$.

It is a well-known result (see, for instance, \cite{KR})
that if a function $u$ of class
$H^{n}_{2}(T)$ is a generalized solution
of \eqref{6.28.1} with initial
condition $u_{0}$, then there exists
a normal solution $\hat{u}$ of \eqref{6.28.1}
of class $H^{n}_{2}(T)$  with initial
condition $u_{0}$ such that $\hat{u}_{t}$
and $u_{t}$
coincide as generalized functions for almost all
$(\omega,t)$.  

This result implies, in particular, that
if a generalized solution of class $H^{n}_{2}(T)$
is such that $u_{t}$ is a continuous 
$H^{m}_{2}$-valued function for some $m$, then
$\hat{u}_{t}=u_{t}$ (a.s.) for all $t\in[0,T]$,
so that $u_{t}$ itself is a continuous 
$H^{n-1}_{2}$-valued function (a.s.)  and thus a normal
solution of class $H^{n}_{2}(T)$.
\end{remark}

Next we need the following technical lemma
which enables us to integrate by parts in $H^{n}_{2}$-spaces.

\begin{lemma}
                                \label{lemma 7.29.1}
Let $n\in \bR$ and 
$\nu\in C^{\infty}_{b}(\bR^{d})$. Then there exists a constant $N$
such that for any $u\in H^{n+1}_{2}$, $k=0,1,...,d_{1}+d_{2}$,
$t\in[0,T]$ and $\omega\in\Omega$ we have
\begin{equation}
                                                       \label{7.30.3}
|(\Lambda^{n}L_{\sigma^{k}_{t}}u,
\Lambda^{n}(\nu u) )_{0}|\leq N\|u\|_{n}^{2},
\end{equation}
\begin{equation}
                                                       \label{7.29.1}
\langle\Lambda^{n}u , \Lambda^{n}(L_{\sigma^{k}_{t}}^{2}u )\rangle
+\|
\Lambda^{n} L_{\sigma^{ k}_{t}}u  \|^{2}_{0}
\leq N\|u\|_{n}^{2},
\end{equation}
where $\langle\cdot,\cdot\rangle$ is the natural pairing
between $H^{1}_{2}$ and $H^{-1}_{2}$.
\end{lemma}

Proof. By obvious reasons we 
may assume that $u\in C^{\infty}_{0}(\bR^{d})$
and we drop the indices $k$ and $t$ to simplify notation. 

We are going to rely on some well-known properties of 
pseudo-differential operators. The order
of a pseudo-differential operator $S$ is a number $n\in\bR$
such that $\Lambda^{-n}S$ and $S\Lambda^{-n}$ are bounded
operators in $\cL_{2}$. If the orders of two operators
$S_{1}$ and $S_{2}$ are $n_{1}$ and $n_{2}$, respectively,
then the order of $[S_{1},S_{2}]=S_{1}S_{2}-S_{2}S_{1}$
is at most $n_{1}+n_{2}-1$. One also knows that
the first order linear differential operators with
 coefficients whose every derivative of any order is bounded
are pseudodifferential operator of order one. 

Observe that if a pseudo-differential operator $S$ is self adjoint,
then for any $u\in C^{\infty}_{0}(\bR^{d})$ 
$$
(L_{\sigma}u,S(\nu u))_{0}=(L_{\sigma}Su,\nu u)_{0}+([S,L_{\sigma}]u,
\nu u)_{0}
$$
$$
=-(Su,\nu L_{\sigma}u)_{0}+(aSu,u)_{0}+([S,L_{\sigma}]u,\nu u)_{0}
$$
$$
=-(S(\nu u),L_{\sigma}u)_{0}-([\nu\cdot,S]u,L_{\sigma}u)_{0}
+(aSu,u)_{0}+([S,L_{\sigma}]u,\nu u)_{0},
$$
where $a$ is a smooth bounded function. It follows that
\begin{equation}
                                                       \label{7.30.1}
(L_{\sigma}u,S(\nu u))_{0}=(1/2)\big[(aSu,u)_{0}+
 ([S,L_{\sigma}]u,\nu u)_{0}- 
 ([\nu\cdot,S]u,L_{\sigma}u)_{0}\big].
\end{equation}
It is important to note for the future that if
the order of $S$ is $2n$, then the order of $[S,L_{\sigma}]$
is at most $2n$, the order of $[\nu\cdot,S]$ is at most $2n-1$, and
consequently
\begin{equation}
                                                       \label{7.30.2}
|(L_{\sigma}u,S(\nu u))_{0}|\leq N\|u\|^{2}_{n}.
\end{equation}
This with $S=\Lambda^{2n}$ yields \eqref{7.30.3}.

By applying \eqref{7.30.1} with  $S=\Lambda^{2n}$ and $\nu\equiv1$,
 we get
$$
(\Lambda^{ n}L_{\sigma}u,\Lambda^{ n}u)_{0}=(1/2)(a\Lambda^{2n}u,u)_{0}+
(1/2)([\Lambda^{2n},L_{\sigma}]u,u)_{0},
$$
which after being polarized yields that
$$
(\Lambda^{ n}L_{\sigma}v,\Lambda^{ n}u)_{0}
+(\Lambda^{ n}L_{\sigma}u,\Lambda^{ n}v)_{0}= (a\Lambda^{2n}u,v)_{0}
$$
$$
+
(a\Lambda^{2n}v,u)_{0}+
 ([\Lambda^{2n},L_{\sigma}]u,v)_{0}+
([\Lambda^{2n},L_{\sigma}]v,u)_{0} 
$$
if $u,v\in C^{\infty}_{0}(\bR^{d})$.
We plug in here $v=L_{\sigma}u$ and obtain
$$
( L^{2}_{\sigma}u,u)_{n}+\|L_{\sigma}u\|^{2}_{n}
=(a\Lambda^{2n}u,L_{\sigma}u)_{0}
$$
$$
+(a\Lambda^{2n}L_{\sigma}u,u)_{0}+
([\Lambda^{2n},L_{\sigma}]u,L_{\sigma}u)_{0}+
([\Lambda^{2n},L_{\sigma}]L_{\sigma}u,u)_{0}.
$$
After introducing the self adjoint operators
$$
S_{1}=a\Lambda^{2n}+(a\Lambda^{2n})^{*},\quad
S_{2}=[\Lambda^{2n},L_{\sigma}]+([\Lambda^{2n},L_{\sigma}])^{*}
$$
we rewrite the last equation as
$$
( L^{2}_{\sigma}u,u)_{n}+\|L_{\sigma}u\|^{2}_{n}
=(L_{\sigma}u,S_{1}u)_{0}+(L_{\sigma}u,S_{2}u)_{0}
$$
and obtain \eqref{7.29.1} owing to \eqref{7.30.2}.
The lemma is proved.

\begin{theorem}
                                     \label{theorem 7.24.1}
In class $H^{n}_{2}(T) $  there exists an (a.s.) unique    
normal  solution $u $
of \eqref{6.28.1} on $[0,T]$ with initial condition $u_{0}$.
Furthermore, there exists a constant $N$ independent of $u,f,g$
such that
\begin{equation}
                                                  \label{7.24.1}
E\sup_{t\leq T}\|u_{t}\|^{2}_{n}\leq
NE\|u_{0}\|^{2}_{n}+NE\int_{0}^{T}\big(\|f_{t}\|^{2}_{n}+
\sum_{k=1}^{d_{1}}\|g^{k}_{t}\|^{2}_{n+1}\big)\,dt.
\end{equation}

\end{theorem}

Proof.  {\em Step 1}.
First we want to
derive an a priori estimate assuming that 
we are given a  normal solution of \eqref{6.28.1}
of class $H^{n+1}(T)$.
We apply the operator $\Lambda^{n}$ to both sides
of \eqref{6.28.1} written in the integral form
 and observe that after that the stochastic
integral will belong to $L_{2}$, whereas the deterministic
integral will belong to $H^{-1}_{2}$. This allows us to apply
It\^o's formula for Banach space valued processes and shows that
$$
d\|u_{t}\|^{2}_{n}=d\|\Lambda^{n}u_{t}\|^{2}_{0}
= I_{t}\,dt
+
2\big(\Lambda^{n}u_{t},
\Lambda^{n}(L_{\sigma^{ k}_{t}}u_{t}+\nu^{k}_{t}u_{t}
+g^{k}_{t})\big)_{0}\,dw^{k}_{t},
$$
where
$$
I_{t}= \langle\Lambda^{n}u_{t},2\Lambda^{n}(L_{t}u_{t}+c_{t}u_{t}
+f_{t})\rangle+\sum_{k=1}^{d_{1}}\|
\Lambda^{n}(L_{\sigma^{ k}_{t}}u_{t}+\nu^{k}_{t}u_{t}
+g^{k}_{t})\|^{2}_{0} 
$$
$$
=I_{t}^{1}+I_{t}^{2}+I_{t}^{3}+2I_{t}^{4}+2I_{t}^{5}+I_{t}^{6},
$$
$$
I_{t}^{1}=\sum_{k=1}^{d_{1}}\big[ 
\langle\Lambda^{n}u_{t}, \Lambda^{n}(L_{\sigma^{k}_{t}}^{2}u
_{t})\rangle
+\|
\Lambda^{n} L_{\sigma^{ k}_{t}}u_{t}\|^{2}_{0}\big],
$$
$$
I_{t}^{2}=(\Lambda^{n}u_{t},2\Lambda^{n} L_{\sigma^{0}_{t} }
u_{t} )_{0},
$$
$$
I_{t}^{3}=(\Lambda^{n}u_{t},2\Lambda^{n}(  c_{t}u_{t}
+f_{t}))_{0},
$$
$$
I_{t}^{4}=  (\Lambda^{n}L_{\sigma^{k}_{t}}u_{t},
\Lambda^{n}(\nu^{k}_{t}u_{t}))_{0},
$$
$$
I_{t}^{5}=  (\Lambda^{n}L_{\sigma^{k}_{t}}u_{t},
\Lambda^{n}g^{k}_{t} )_{0},
$$
$$
I_{t}^{6}=
\sum_{k=1}^{d_{1}}\|
\Lambda^{n} (\nu^{k}_{t}u_{t}
+g^{k}_{t})\|^{2}_{0}.
$$

The term $I^{1}_{t}$ is estimated in \eqref{7.29.1}  and
  $I^{2}_{t}$ in \eqref{7.30.3} (with $\nu\equiv1$),
which also provides an estimate for $I^{4}_{t}$. Almost obviously
$$
|I^{3}_{t}|+|I^{6}_{t}|\leq N\|u_{t}\|^{2}_{n}+N\|f_{t}\|^{2}_{n}
+N\sum_{k=1}^{d_{1}}\|g^{k}_{t}\|^{2}_{n},
$$
where and below we denote by $N$ various
constants independent of $u,f,g^{k}$, $t$, and $\omega$.
Finally,
$$
I_{t}^{5}=  (L_{\sigma^{k}_{t}}\Lambda^{n}u_{t},
\Lambda^{n}g^{k}_{t} )_{0}+([\Lambda^{n},L_{\sigma^{k}_{t}}]u_{t},
\Lambda^{n}g^{k}_{t})_{0}
$$
$$
= (\Lambda^{n}u_{t},(L_{\sigma^{k}_{t}})^{*}
\Lambda^{n}g^{k}_{t} )_{0}+([\Lambda^{n},L_{\sigma^{k}_{t}}]u_{t},
\Lambda^{n}g^{k}_{t})_{0}
$$
and, since the order of the operator $[\Lambda^{n},L_{\sigma^{k}_{t}}]$
is at most $n$, we have
$$
|I^{5}_{t}|\leq  N\sum_{k=1}^{d_{1}}\|g^{k}_{t}\|^{2}_{n+1}.
$$

Upon collecting our estimates we conclude that
$$
d\|u_{t}\|^{2}_{n}\leq N\big(
\|u_{t}\|^{2}_{n}+\|f_{t}\|^{2}_{n}
+ \sum_{k=1}^{d_{1}}\|g^{k}_{t}\|^{2}_{n+1}\big)\,dt
$$
\begin{equation}
                                                       \label{7.30.4}
+
2\big(\Lambda^{n}u_{t},
\Lambda^{n}(L_{\sigma^{ k}_{t}}u_{t}+\nu^{k}_{t}u_{t}
+g^{k}_{t})\big)_{0}\,dw^{k}_{t}.
\end{equation}

{\em Step 2. Uniqueness}.
Now assume that we are given two natural 
solutions of \eqref{6.28.1} of class $H^{n}_{2}(T)$  
with the same initial condition.
Then for their difference, say $u_{t}$ we have
$$
d\|u_{t}\|^{2}_{n-1}\leq N_{1} 
\|u_{t}\|^{2}_{n-1} \,dt+dm_{t},
$$
where $m_{t}$ is a local martingale. Next, comparing the
differentials we obtain
$$
\|u_{t}\|^{2}_{n-1}e^{-N_{1}t}\leq \int_{0}^{t}e^{-N_{1}s}\,dm_{s}.
$$
Since the right-hand side is a local martingale starting at zero
and the left-hand side is nonnegative, it follows, that
the right-hand side is zero as is the left-hand side, which
proves uniqueness.

{\em Step 3}. Here we prove   \eqref{7.24.1} as an a priori
estimate under the assumptions of Step 1. 
We follow by now an absolutely standard and well-known way.
 With $N$ from \eqref{7.30.4}
we have
$$
d(e^{-Nt}\|u_{t}\|^{2}_{n})\leq N\big(
 \|f_{t}\|^{2}_{n}
+ \sum_{k=1}^{d_{1}}\|g^{k}_{t}\|^{2}_{n+1}\big)\,dt+dm_{t},
$$
where $m_{t}$ is a local martingale. Since the left-hand side
is nonnegative, for any $t\in[0,T]$,
$$
e^{-Nt}E\|u_{t}\|^{2}_{n}\leq E\|u_{0}\|^{2}_{n}+ NE\int_{0}^{t}\big(
 \|f_{s}\|^{2}_{n}
+ \sum_{k=1}^{d_{1}}\|g^{k}_{s}\|^{2}_{n+1}\big)\,ds,
$$
\begin{equation}
                                                       \label{7.30.5}
 \sup_{t\in[0,T]}
E\|u_{t}\|^{2}_{n}\leq NE\|u_{0}\|^{2}_{n}+NE\int_{0}^{T}\big(
 \|f_{s}\|^{2}_{n}
+ \sum_{k=1}^{d_{1}}\|g^{k}_{s}\|^{2}_{n+1}\big)\,ds.
\end{equation}

Next, notice that since $\|u_{t}\|_{n}$ is continuous
$$
\tau_{m}:=T\wedge\inf\{t\geq0:\|u_{t}\|_{n}\geq m\}
$$
are stopping times and $\tau_{m}\uparrow T$ as $m\to\infty$.
By  Davis's inequality \eqref{7.30.4}
and \eqref{7.30.5}
imply that
$$
E\sup_{t\leq\tau_{m}}\|u_{t}\|^{2}_{n}
\leq NE\|u_{0}\|^{2}_{n}+NE\int_{0}^{T}\big(
 \|f_{s}\|^{2}_{n}
+ \sum_{k=1}^{d_{1}}\|g^{k}_{s}\|^{2}_{n+1}\big)\,ds
$$
\begin{equation}
                                                       \label{7.30.6}
+6E\big(\int_{0}^{\tau_{m}}\sum_{k=1}^{d_{1}}\big|
\big(\Lambda^{n}u_{t},
\Lambda^{n}(L_{\sigma^{ k}_{t}}u_{t}+\nu^{k}_{t}u_{t}
+g^{k}_{t})\big)_{0}\big|^{2}\,dt\big)^{1/2}.
\end{equation}
By the above what is inside the square by magnitude
is dominated by 
$$
N\|u_{t}\| _{n}(\|u_{t}\| _{n}+\|g^{k}_{t}\| _{n}).
$$ 
Hence the last term in \eqref{7.30.6} is less than
$$
NE\big(\int_{0}^{\tau_{m}}
\|u_{t}\|^{2}_{n}\big(\|u_{t}\|^{2}_{n}+
\sum_{k=1}^{d_{1}}\|g^{k}_{t}\|^{2}_{n}\big)\,dt\big)^{1/2}
$$
$$
\leq NE\big(\sup_{t\leq\tau_{m}}\|u_{t}\|_{n}
 \int_{0}^{\tau_{m}}
 \big(\|u_{t}\|^{2}_{n}+
\sum_{k=1}^{d_{1}}\|g^{k}_{t}\|^{2}_{n}\big)\,dt\big)^{1/2}
$$
$$
\leq(1/2)E\sup_{t\leq\tau_{m}}\|u_{t}\|^{2}_{n}
+NE\int_{0}^{\tau_{m}}
 \big(\|u_{t}\|^{2}_{n}+
\sum_{k=1}^{d_{1}}\|g^{k}_{t}\|^{2}_{n}\big)\,dt,
$$
which after coming back to \eqref{7.30.6}, using again \eqref{7.30.5},
and sending $m\to\infty$, by Fatou's lemma yields the a priori estimate
\eqref{7.24.1}.

{\em Step 4. Existence in a particular case}.
If the norms on the right in \eqref{7.24.1}
are finite, $u_{0}\in L_{2}(\Omega,\cF_{0},H^{n+1}_{2})$,
 and our equation is uniformly nondegenerate,
then (see, for instance, \cite{Kr99}) there exists a unique
normal  solution of our problem of class
$H^{n+2}_{2}(T)$. For this solution \eqref{7.24.1} is valid.

If the norms on the right in \eqref{7.24.1}
are finite, but $u_{0}\in L_{2}(\Omega,\cF_{0},H^{n }_{2})$,
and there is no nondegeneracy assumption,
  we approximate $u_{0}$ in the 
$ L_{2}(\Omega,\cF_{0},H^{n }_{2})$-norm by
a sequence $u^{m}_{0}\in L_{2}(\Omega,\cF_{0},H^{n+1}_{2})$,
$m=1,2,...$, and add into the right-hand side of
\eqref{6.28.1} the term $(1/m)\Delta u_{t}\,dt$ to make the equation
uniformly nondegenerate.
Denote by $u^{m}_{t}$ the solutions of the so modified problems.
Then \eqref{7.24.1} will hold with $N$ independent of $m$
because no constant of nondegeneracy was involved in the
derivation of \eqref{7.24.1}.

According to \eqref{7.24.1} the sequence $u^{m}_{t}$ is bounded
in   $\bH^{n}_{2}(T)$. In particular, $(1/m)\Delta u_{t}\to0$
in $\bH^{n-2}_{2}(T)$. Having in mind his fact and
applying \eqref{7.24.1} to the difference $u^{m}_{t}-u^{k}_{t}$
and $n-2$ in place of $n$ we see
 that the sequence $u^{m}_{t}$  is Cauchy
  in   the space with norm, whose square is given by
$$
E\sup_{t\leq T}\|u_{t}\|_{n-2}^{2},
$$
in particular, in
$\bH^{n-2}_{2}(T)$. Let $u_{t}$ be its limit in
$\bH^{n-2}_{2}(T)$ such that
\begin{equation}
                                              \label{7.30.06}
E\sup_{t\leq T}\|u_{t}-u^{m}_{t}\|_{n-2}^{2}\to0
\end{equation}
as $n\to\infty$. 
 Then equation \eqref{6.28.1}
in the integral form holds in $\bH^{n-4}_{2}(T)$.
Now, since the sequence $u^{m}_{t}$ is bounded
in   $\bH^{n}_{2}(T)$ and converges to $u_{t}$ in $\bH^{n-2}_{2}(T)$,
  $u_{t}\in \bH^{n}_{2}(T)$. After that we apply a classical result
saying that if $u_{t}\in \bH^{n}_{2}(T)$ satisfies \eqref{7.24.1}
in  generalized sense
with initial condition $u_{0}\in H^{n-1}_{2}$, then (a.s.)
$u_{t}$ is a continuous $H^{n-1}_{2}$-valued function
(see Remark \ref{remark 7.31.2}).

To establish \eqref{7.24.1} for thus found normal solution
take a sequence $\phi^{r}\in C^{\infty}_{0}(\bR^{d})$
such that it is dense in the unit ball of $H^{ n}_{2}$. Then
owing to \eqref{7.24.1} write for any $j=1,2,...$
\begin{equation}
                                                      \label{7.30.7}
E\sup_{t\leq T}\max_{r=1,...,j}
(\Lambda^{n} u^{m}_{t},\Lambda^{n} \phi^{r})^{2}_{0} \leq
NE\|u^{m}_{0}\|^{2}_{n}+I,
\end{equation}
where $I$ is the second term on the right in \eqref{7.24.1}.
Since 
$$
(\Lambda^{n} (u^{m}_{t}-u_{t}),\Lambda^{n} \phi^{r})_{0}^{2}
\leq \|\Lambda^{n-2}(u^{m}_{t}-u_{t})\|^{2}_{0}
\|\Lambda^{n+2} \phi^{r})_{0}^{2},
$$
estimate \eqref{7.30.06} allows us to conclude from
\eqref{7.30.7} that
$$
E\sup_{t\leq T}\max_{r=1,...,j}
(\Lambda^{n} u _{t},\Lambda^{n} \phi^{r})^{2}_{0} \leq
NE\|u _{0}\|^{2}_{n}+I.
$$
By letting $j\to\infty$ and using the monotone convergence
theorem and the fact that $\phi^{r}$ are 
  dense in the unit ball of $H^{ n}_{2}$, we get \eqref{7.24.1}.

{\em Step 5. Existence in the general case}.
The first assertion of the theorem
in the general case is proved as always by using stopping times
like
$$
\gamma_{m}=\inf\{t\geq0:
\|u_{0}\|^{2}_{n}+ \int_{0}^{t}\big(\|f_{s}\|^{2}_{n}+
\sum_{k=1}^{d_{1}}\|g^{k}_{s}\|^{2}_{n+1}\big)\,ds\geq m\}.
$$
The theorem is proved.

In the remaining part of this section by $u_{t}$
we mean the normal solution from 
Theorem 
\ref{theorem 7.24.1}. We remind the reader that the common way of saying that
a generalized function in a domain is smooth
means that there is a smooth function which, 
as
a as generalized function, coincides
with the given generalized one in the domain
under consideration.

Theorem \ref{theorem 7.24.1}
and Sobolev embedding theorems immediately
imply the following.

\begin{corollary}
                                         \label{corollary 8.1.1}
Suppose that Assumption \ref{assumption 6.28.1}  
is satisfied for all $n$. Then (a.s.) the solution $u_{t}$
 is infinitely differentiable with respect to $x$
and every its derivative is a bounded continuous
function of $(t,x)$.
\end{corollary}

\begin{corollary}
                                         \label{corollary 8.1.2}
Suppose that Assumption \ref{assumption 6.28.1}  
is satisfied for all $n$. Let
$D$ be a domain in $\bR^{d}$ with $\partial D\ne\emptyset$
and assume that
for $x\in D$, $t\in[0,T]$, and $\omega\in\Omega$ we have
$$
u_{0}(x) ,c_{t}(x),f_{t}(x)\leq 0,\quad \nu^{k}_{t}(x)=g^{k}_{t}(x)=0,
\quad k=1,...,d.
$$
Then (a.s.) for all $t\in[0,T]$ on $D$ we have
\begin{equation}
                                                   \label{8.1.2}
u_{t}\leq\max_{s\leq t}\max_{\partial D}(u_{s})_{+}.
\end{equation}
\end{corollary}

This result follows from Theorem 1.2 of \cite{Kr07}
in which one takes $\xi=0$, $\bar{u}\equiv 1$,
$\bar{f} \equiv -c$,
$\bar{f}^{i}\equiv0$, $\bar{g}\equiv0$, and $\rho_{t}$
equal the right-hand side of \eqref{8.1.2} plus a constant
$\varepsilon>0$. One adds $\varepsilon$ to be sure that
$$
(u_{t}-\rho_{t}\bar{u}_{t})_{+}=
(u_{t}-\rho_{t} )_{+}
$$
vanishes near the boundary of $D$ and hence belongs to
$\WO^{1}_{2}(D)$. Then after applying
Theorem 1.2 of \cite{Kr07} one sets $\varepsilon\downarrow0$.

\begin{theorem}
                                            \label{theorem 8.1.1}
Take an $R\in[0,\infty)$
and suppose that $\sigma^{k}_{t}(x) =0$ and $\nu^{k}_{t}(x) =0$ 
for
 $k=1,...,d_{1}$, $t\in[0,T]$, and $\omega\in\Omega$
as long as $|x|>R$.
Also assume that $g^{k}\equiv0$ for
 $k=1,...,d_{1}$. Then  (a.s.)
there exists a (random finite) 
constant $N$ independent of $f$ and $u_{0}$ such that
\begin{equation}
                                                  \label{8.1.4}
 \sup_{t\leq T}\|u_{t}\|^{2}_{n}\leq
N \|u_{0}\|^{2}_{n}+N \int_{0}^{T} \|f_{t}\|^{2}_{n} \,dt.
\end{equation}

\end{theorem}

Proof. For smooth $\bR^{d}$-valued functions $\sigma(x)$ on
$\bR^{d}$ (whose points are always considered as column vectors)
by 
$D\sigma$ we mean a matrix with entries $
(D\sigma)^{ij}=D_{j}\sigma^{i}$ and if we are given two such
functions $\sigma$ and $\gamma$, then we set
\begin{equation}
                                             \label{8.2.6}  
D\sigma\gamma:=[D\sigma]\gamma.
\end{equation}
 Consider the equation
\begin{equation}
                                             \label{6.28.4}  
x_{t}=x-\int_{0}^{t}\sigma^{k}_{s}( x_{s})\,dw^{k}_{s}
-\int_{0}^{t}b_{t}(x_{s})\,ds,
\end{equation}
where
$$
b_{t}(x)=\sigma^{0}_{t}(x)-(1/2)\sum_{k=1}^{d_{1}}D\sigma^{k}
_{t}(x)\sigma^{k}_{t}(x).
$$
As  is well known (see, for instance,  \cite{Ku90} for more  
advanced treatment of the subject or see \cite{Kr13_1}),
 there exists
a function $X_{t}(x)$ on $\Omega\times[0,T]\times\bR^{d}$,
such that

(i) it is continuous in $(t,x)$ for any $\omega$
along with each derivative of $X_{t}(x)$ of any order
with respect to $x$,

(ii) for each $\omega$ and $t$ the mapping $x\to X_{t}(x)$
of $\bR^{d}$ to $\bR^{d}$
is one-to-one and onto and its inverse mapping $X^{-1}_{t}(x)$
has bounded and continuous in $(t,x)$ derivatives
of any order with respect to $x$ for any $\omega$.

(iii) it is  $\cF_{t}$-adapted for any $x$,

(iv) for each $x$ with probability one it satisfies
\eqref{6.28.4} for all $t\in[0,T]$,  

Observe that $X_{t}(x)=x$ for $|x|\geq R$, and $X^{-1}_{t}(x)
=x$ for all $t\in[0,T]$ if $|x|$ is large enough (depending on
$\omega$).

Next, define the operations  ``hat'' and ``check''
which transform  any
  function
$\phi_{t}(x)$ into
$$
\hat{\phi}_{t}(x):=\phi_{t}( X_{t}(x)),\quad \check{\phi}
=\phi_{t}( X^{-1}_{t}(x)).
$$
Also define $\rho_{t}(x)$ from the equation
$$
\rho_{t}(X_{t}(y))\det DX_{t}(y)=1
$$
and observe that by the change of variables formula
\begin{equation}
                                                  \label{6.30.1}
\int_{\bR^{d}}F(X_{t}(y))\phi(y)\,dy=
\int_{\bR^{d}}F(x)\check{\phi}_{t}(x)\rho_{t}(x)
\,dx,
\end{equation}
 whenever at least one side of the equation
makes sense. Finally, define the mapping   ``bar'' which transforms    
any $\bR^{d}$-valued function $\sigma_{t}(x)$ into
\begin{equation}
                                             \label{7.23.1}
\bar{\sigma}_{t}(x)=Y_{t}(x)\hat{\sigma}_{t}(x)
=Y_{t}(x)\sigma_{t}( X_{t}(x)),
\end{equation}
where
$$
Y=(DX)^{-1}.
$$

By Corollary 6.5 of \cite{Kr13_1} (also see Remark
\ref{remark 7.31.2}) the function $\hat{u}_{t}$
is well defined and is a normal solution of class
$H^{n}_{2}(T)$ of the equation
$$
 d\hat{u}_{t} = \big[\sum_{k=1}^{d_{2}}
 L_{\bar \sigma_{t}^{d_{1}+k}}^{2} \hat{u}_{t}
+\hat{c}_{t}\hat{u}_{t}+\hat{f}_{t}   \big]\,dt
+ \hat{u}_{t}\hat{\nu}^{k}_{t}  \,dw^{k}_{t}.
$$

By using Kolmogorov's continuity theorem for random fields
one easily shows that there exists a function 
$I_{t}(x)=I_{t}(\omega,x)$ which along with each its derivative
of any order
with respect to $x$ is
continuous with respect to $(t,x)\in[0,T]\times\bR^{d}$
for each $\omega$ and such that for each $(t,x)$ 
with probability one
$$
I_{t}(x)=\int_{0}^{t}\hat{\nu}^{k}_{s}(x)  \,dw^{k}_{s}.
$$ 
Then define
$$
\gamma_{t}(x)=\exp\big[-I_{t}(x)-(1/2)\sum_{k=1}^{d_{1}}\int_{0}^{t}
|\hat{\nu}^{k}_{s}(x) |^{2}\,ds\big],\quad
v_{t} =\hat{u}_{t}\gamma_{t}.
$$
Of course, $\gamma_{t}(x)=1$ if $|x|\geq R$.
We want to apply It\^o's formula to write an equation
for $v_{t}$, that is, for any $\phi\in C^{\infty}_{0}(\bR^{d})$
 find the stochastic differential
of
$$
I_{t} (\phi)
:=(v_{t},\phi)=(\hat{u}_{t},\gamma_{t}\phi)
=\big(\Lambda^{-(n-1)}\hat{u}_{t},\Lambda^{ n-1 }(\gamma_{t}\phi)\big
)_{0}.
$$
Here   $\xi_{t}:=\Lambda^{-(n-1)}\hat{u}_{t}$
and $\eta_{t}:=\Lambda^{ n-1 }(\gamma_{t}\phi)$ are continuous
$\cL_{2}$-valued processes, admitting stochastic differentials
such that the classic formula for the squared norm is
applicable. Then this formula is also applicable to
$\xi_{t}+\lambda\eta_{t}$ for any number $\lambda$.
By comparing the coefficients of $\lambda$
in $\|\xi_{t}+\lambda\eta_{t}\|^{2}_{0}$ and in the formula
we obtain the stochastic differential of $(\xi_{t},\eta_{t})_{0}$
that is of $I_{t} (\phi)$.

In this way we get that with probability one
for all $t\in[0,T]$
\begin{equation}
                                                     \label{8.2.5}
 v_{t}=v_{0}+\int_{0}^{t}
\big[\sum_{k=1}^{d_{2}}
\gamma_{s} L_{\bar \sigma_{s}^{d_{1}+k}}^{2}(\gamma_{s}^{-1}v_{s})
+\hat{c}_{s}v_{s}+\gamma_{s}\hat{f}_{s}   \big]\,ds.
\end{equation}
Fix an $\omega$ such that \eqref{8.2.5} holds
for all $t\in[0,T]$, 
$$
\int_{0}^{T}\|u_{t}\|^{2}_{n}\,dt<\infty,
$$
 $u_{0}(\omega,\cdot)
\in H^{n}_{2}$, and
$u_{t}(\omega,\cdot)$ is an $H^{n-1}_{2}$-valued
continuous function. Then \eqref{8.2.5} becomes a deterministic
equation, to which Theorem \ref{theorem 7.24.1} is applicable
because the differential operators in \eqref{8.2.5}
can be rewritten in a canonical form as in
\eqref{6.28.1} with coefficients satisfying 
Assumption \ref{assumption 6.28.1}. Then applying 
Theorem \ref{theorem 7.24.1} to each particular $\omega$,
the set of which has full probability, we obtain
\eqref{8.1.4}. The theorem is proved.

\begin{lemma}
                                               \label{lemma 8.19.1}
Take
$$
f_{t}(x)=u_{0}(x)=(1+|x|^{2})^{-d},\quad g^{k}_{t}(x)=0
$$
and call $v_{t}$  the normal solution
of \eqref{6.28.1} on $[0,T]$ with so prescribed data.
By Corollary \ref{corollary 8.1.1} the function $v_{t}$
is (a.s.) infinitely differentiable with respect to $x$
and every its derivative is a bounded continuous
function of $(t,x)$. We assert that with probability one
for every $r\in(0,\infty)$ there exists a (random constant)
$\varepsilon>0$ such that $v_{t}(x)\geq\varepsilon$ for
$t\in[0,T]$ and $x\in \bar B_{r}$.
\end{lemma}

Proof. 
First observe that $v_{t}\geq0$ by the maximum principle.
We take an $R>r$ and concentrate on  equation
\eqref{6.28.1}  only for $x\in B_{r}$. Then equation
in $B_{r}$ will still hold if we cut off all $\sigma^{k}_{t}$
and  $\nu^{k}_{t}$
outside $B_{r}$ so that they will vanish outside $B_{R}$.
For simplicity of notation we  assume that
$\sigma^{k} =0$ and $\nu^{k} =0$ outside $B_{R}$ for
 $k=1,...,d_{1}$ for the original coefficients. 

Then making the same transformations as in the proof of Theorem
\ref{theorem 8.1.1} we come to the conclusion that 
for almost any $\omega$

(i) equation
\eqref{8.2.5} holds on $ \{(t,x):t\in [0,T],x\in X^{-1}_{t}(
\bar B_{R})\}$
with $\hat{v}_{t}\gamma_{t}$ in place of $v_{t}$. 

For each 
$\omega$ this is a deterministic (degenerate) parabolic equation with
bounded coefficients. Furthermore, for almost any $\omega$

(ii) $\gamma_{s}\hat{f}_{s}>0$ and $\hat{v}_{t}\gamma_{t}\geq0$
on $ \{(t,x):t\in [0,T],x\in X^{-1}_{t}(\bar B_{R})\}$.

Now we want to use the maximum principle to show that
$v_{t}(x)$ cannot take zero value in $[0,T]\times\bar B_{r}$
whenever $\omega$ is such that (i) and (ii) hold. Were the coefficients
of \eqref{8.2.5} continuous in $t$, this would be just a trivial matter.
In our case we still need a little argument. Assume the contrary:
there is a point $(t_{0},x_{0})\in [0,T]\times\bar B_{r}$
such that $v_{t_{0}}(x_{0})=0$. Obviously, $t_{0}>0$.
Then, since $v_{t}\geq0$, for $t\in[0,t_{0})$ we get
$$
0=\hat{v}_{t}(x_{0})\gamma_{t}(x_{0})+
\int_{t}^{t_{0}}
\big[\sum_{k=1}^{d_{2}}
\gamma_{s} L_{\bar \sigma_{s}^{d_{1}+k}}^{2}  \hat{v}_{s} 
+\hat{c}_{s}v_{s}+\gamma_{s}\hat{f}_{s}   \big](x_{0})\,ds
$$
$$
\geq \int_{t}^{t_{0}}
\big[\sum_{k=1}^{d_{2}}
\gamma_{s} L_{\bar \sigma_{s}^{d_{1}+k}}^{2}  \hat{v}_{t_{0}} 
+\hat{c}_{s}v_{t_{0}}+\gamma_{s}\hat{f}_{s}   \big](x_{0})\,ds
+I_{t},
$$
where $I_{t}$ is defined as the difference of the above
two integrals. Since $v_{t}(x)$ and its 
 derivatives with respect to $x$  
are continuous with respect to $t$, we have
$$
\lim_{t\uparrow t_{0}}\frac{1}{t_{0}-t}I_{t}=0.
$$
By taking into account that $v_{t_{0}}(x_{0})=0$,
the first order derivatives of $v_{t_{0}}(x)$ vanish at $x_{0}$,
and the Hessian is nonnegative at $x_{0}$ we conclude that
$$
0\geq \nlimsup_{t\uparrow t_{0}}\frac{1}{t_{0}-t}
\int_{t}^{t_{0}}\gamma_{s}\hat{f}_{s}(x_{0})\,ds,
$$
which is impossible because $\gamma_{s}\hat{f}_{s}(x_{0})$
is strictly positive on $[0,t_{0}]$.
The lemma is proved.

\begin{theorem}
                                         \label{theorem 8.13.1}
Suppose that Assumption \ref{assumption 6.28.1}  
is satisfied for all $n$. 
Also assume that there is an $r\in(0,\infty)$ such that $
u_{0}(x) ,f_{t}(x)\leq 0$, and 
$g^{k}_{t}(x)= 0$ for
 $k=1,...,d_{1}$ and $x\not\in B_{r}$. Then  (a.s.)
there exists a (random finite) 
constant $N$ independent of $f,g^{k}$, and $u_{0}$ such that
 (a.s) for $t\in[0,T]$ and $|x|\geq r$ we have
\begin{equation}
                                                   \label{8.13.1}
u_{t}(x)\leq N\max_{s\leq t}\max_{\partial B_{r}}(u_{s})_{+}.
\end{equation}
\end{theorem}

Proof. Take $v_{t}$ from Lemma \ref{lemma 8.19.1} and set
$$
\rho_{t}=\max_{s\leq t}\max_{|x|=r}(u_{t}(x)/v_{t}(x))_{+}.
$$
By  Lemma \ref{lemma 8.19.1} the process $\rho_{t}$
is finite, nonnegative, increasing, and continuous with probability one. Furthermore,
$(u_{t}-\rho_{t}v_{t} )_{+}$ vanishes on $\partial B_{r}$.
This along with the fact that $
u_{0}(x) =f_{t}(x)\leq 0$ and 
$g^{k}_{t}(x)= 0$ for
 $k=1,...,d_{1}$ and $x\not\in B_{r}$  by Theorem 1.2 of \cite{Kr07}
implies that $ u_{t}-\rho_{t}v_{t}\leq 0$ in $[0,T]\times 
B^{c}_{R}$ (a.s.),
which obviously proves the theorem.

\mysection{Hypoellipticity}
                                          \label{section 8.26.1}

Recall the notation associated with \eqref{8.2.6}
and for two smooth $\bR^{d}$-valued functions $\sigma$
and $\gamma$ on $\bR^{d}$ set, as usual,
$$
[\sigma,\gamma]=D\gamma\sigma-D\sigma\gamma.
$$

\begin{assumption}
                                      \label{assumption 8.2.2}
Assumption \ref{assumption 6.28.1} (i)  
is satisfied,
 Assumption \ref{assumption 6.28.1} (ii)  
is satisfied for all $n$  and, for an $n$,
the function $u_{0}$ is an $\cF_{0}$-measurable $H^{n}_{2}$-valued
function on $\Omega$.
\end{assumption}

Fix an $R_{0}\in(0,\infty)$ and
 introduce collections of $\bR^{d}$-valued functions
defined on $\Omega\times[0,T]\times B_{R_{0}}$ inductively as
$\bL_{0}=\{\sigma^{d_{1}+1},...,\sigma^{{d_{1}+d_{2}}}\}$,  
$$
\bL_{n+1}=\cL_{n}\cup\{[\sigma^{d_{1}+k},M]:k=1,...,d_{2},
M\in\bL_{n}\},\quad n\geq0.
$$ 

For any multi-index $\alpha=(\alpha_{1},...,\alpha_{d})$,
$\alpha_{i}\in\{0,1,...\}$, introduce as usual
$$
D^{\alpha}=D^{\alpha_{1}}_{1}\cdot...\cdot D^{\alpha_{d}}_{d},
\quad |\alpha|=\alpha_{1}+...+\alpha_{d}.
$$
Also define
 $BC^{\infty}_{b}$ as the set of real-valued  
measurable
 functions $a$ on $\Omega\times[0,T]\times\bR^{d}$
such that for each $t\in[0,T]$ and $\omega\in\Omega$,
 $a_{t}(x)$ is infinitely differentiable with respect to $x$, and
  for any $\omega\in\Omega$ and multi-index $\alpha$ we have
$$
\sup_{ t,x\in [0,T]\times\bR^{d}}|D^{\alpha}a_{t}( x)|<\infty.
$$

  Finally we denote by
$\Lie_{n}$ the set of (finite) linear combinations
of elements of $\bL_{n}$ with  
coefficients which  are
of class $BC^{\infty}_{b}$.
Observe that the vector-field $\sigma^{0}$ is {\em not\/}
explicitly included
into $\Lie_{n}$. Finally, fix an $S\in[0,T)$
and introduce
 $$
G=(S,T)\times B_{R_{0}}.
$$

\begin{assumption}
                                      \label{assumption 8.2.1}
For every $\omega\in\Omega$  
  and $\zeta\in C^{\infty}_{0}(B_{R_{0}})$
there exists an $n\in\bR$ such that
 we have  $\xi I_{[S,T]} \zeta \in\Lie_{n}$ 
for any $\xi\in\bR^{d}$. 
\end{assumption}

The following result will be used quite often.
It is  a particular case of Theorem 2.3 of \cite{Kr13_1}. 
By $u$ in this theorem 
and everywhere below in this section we mean the normal solution
which exists due to Theorem \ref{theorem 7.24.1}.

\begin{theorem}
                                 \label{theorem 7.22.1} 

Take $ s_{0}\in(S,T)$,
 $r\in(0,R_{0})$ and   take a $\zeta\in C^{\infty}_{0}(B_{R_{0}})$
such that $\zeta=1$ on a neighborhood of
 $ \bar B_{r}$.  Then

(i) with probability one $u_{t}( x)$ is infinitely differentiable 
with respect to $x$ for
$(t,x)\in (S,T]\times B_{R_{0}}$ and each derivative is 
a continuous function
in $(S,T]\times B_{R_{0}}$.

(ii)  
for any multi-index $\alpha$ and $l$ such that
\begin{equation}
                                                 \label{7.16.1} 
2(l-|\alpha|-2)>d+1,
\end{equation}
with probability one
there exists a (random, finite)
constant $N$, independent of $u,f$, and $g^{k}$, such that
\begin{equation}
                                                 \label{7.22.2}
\sup_{(t,x)\in[s_{0},T]\times B_{r}}
|D^{\alpha}u_{t}(x)|^{2}\leq N
\int_{S}^{T}\big[\|f_{t}\zeta \|_{H^{l}_{2}}^{2}
+\| u_{t}  
\zeta \|_{H^{n}_{2} }^{2}\big]\,dt,
\end{equation}
provided that $g^{k}_{t}\zeta\equiv0$, $k=1,...,d_{1}$.

If we additionally assume that $u_{S} $ is infinitely differentiable
in $B_{R_{0}}$ for every $\omega$,
then assertion (i) holds true with $[S,T]\times B_{R_{0}}$
in place of $(S,T]\times B_{R_{0}}$, and  assertion (ii) with
$s_{0}=S$ 
 if we add to the right-hand side
of \eqref{7.22.2} a constant (independent of $u$) times
$\|\zeta u_{S}\|^{2}_{H^{l+1}_{2}}$.

\end{theorem}

Here is a generalization of the corresponding results of paper
\cite{CaM}, where there is no stochastic terms in the equation.
This is a generalization because no continuity hypothesis
in time on the coefficients is imposed.
The types of H\"ormander's condition imposed in \cite{CaM} and here
coincide. It is worth noting, however, that
 the result of \cite{CaM}
bears on the equation formally adjoint to the one we consider
when there is no stochastic terms and no dependence
on $\omega$. Such equations have the same form
as ours and have the same $\Lie_{n}$.
Recall that $\cHO^{n}_{2}(B_{r})$
is introduced before Assumption \ref{assumption 6.28.1}.

\begin{theorem}
                                        \label{theorem 8.18.1}
Suppose that $S=0$,  
 take $s_{0}\in(0,T)$, $r \in(0,R_{0})$,
and suppose that $u_{0}\in \cHO^{n}_{2}(B_{r})$.
 Then $u$ is
(a.s.)  infinitely differentiable 
with respect to $x$ for
$(t,x)\in ([s_{0},T]\times \bR^{d})
\cup([0,T]\times B_{R_{0}}^{c})$ and each 
derivative is a bounded continuous function
in $([s_{0},T]\times \bR^{d})\cup([0,T]\times B_{R_{0}}^{c})$. 
\end{theorem}

Proof. Take $r<r_{1}<r_{2}<R$ and observe that
$B_{r_{2}}\setminus B_{r_{1}}$ can be covered by a finite number of balls
lying inside $B_{R_{0}}\setminus B_{r}$, where $u_{0}=0$.
By applying Theorem \ref{theorem 7.22.1}   to each such ball
we conclude that $u_{t}$ is infinitely differentiable with respect 
$x$ in  $[0,T]\times(B_{r_{2}}\setminus B_{r_{1}})$
and each its derivatives is bounded and continuous in
$[0,T]\times(B_{r_{2}}\setminus B_{r_{1}})$.

Take a $\zeta\in C^{\infty}_{0}(B_{r_{2}})$
such that $\zeta=1$ in a neighborhood of $\bar B_{r_{1}}$
and set $\eta=1-\zeta$. Then $u_{t}\eta$ satisfies an equation
similar to \eqref{6.28.1} but with different $f$ and $g^{k}$
which are obtained by adding to the original ones
$u_{t}$ or its first-order derivatives
multiplied by $C^{\infty}_{0}(\bR^{d})$ functions
which vanish outside $B_{r_{2}}\setminus B_{r_{1}}$.
The initial condition for $u_{t}\eta$ is obviously zero.
By Theorem \ref{theorem 7.24.1} and embedding theorems 
we conclude that
(a.s.) $u_{t}\eta$ is infinitely differentiable 
with respect to $x$ for
$(t,x)\in [0,T]\times \bR^{d}$ and each 
derivative is a bounded continuous function
in $[0,T]\times \bR^{d}$.

The function $u_{t}\zeta$ satisfies an equation
with the properties similar to those of the equation
for $u_{t}\eta$ and by Theorem \ref{theorem 7.22.1}  (a.s.)
is
infinitely differentiable 
with respect to $x$ for
$(t,x)\in [s_{0},T]\times B_{R_{0}}$ and each 
derivative is a continuous function
in $[s_{0},T]\times B_{R_{0}}$. This proves the present theorem
since $u_{t}=u_{t}\eta+u_{t}\zeta$ and $\zeta=0$ outside $B_{R_{0}}$.

\mysection{Applications to filtering problems}
                                         \label{section 9.10.1}

Here we come back to system \eqref{eq3.2.14}.
 Let $K,\delta>0$ be some fixed constants. We denote by $\theta$
and $\Theta$ the matrix-valued functions having $\theta^{k}$ and
$\Theta^{k}$, respectively, as their $k$th columns.

\begin{assumption}
                                             \label{asm3.2.15}
The functions $b$, $\theta $, $B$, and $\Theta $ are
 Borel measurable and bounded
functions of their arguments. Each of them satisfies 
the Lipschitz
condition in $z$ with constant $K $. These functions
are infinitely differentiable with respect to $x$
and each derivative of any order is a bounded function
of $(t,z)$.
\end{assumption}

\begin{assumption}
                                             \label{assumption 8.21.1}
For all $\xi\in \bR^{d'}$ and $(t,y)$
$$
 |\Theta^{ *}(t,y)\xi| \geq
\delta|\xi|.
$$
\end{assumption}

Notice that in \cite{CM} and \cite{Ch} this condition is satisfied
since there $\Theta=(0,I)$, where $I$ is the unit
$d'\times d'$-matrix.
In \cite{KS} there is no such condition because
of a very peculiar filtering problem  considered when the coefficients
of the equation defining $y_{t}$ are independent of $x_{t}$.

\begin{remark}
Owing to Assumption \ref{assumption 8.21.1}, $d'\leq d_{1}$,
the symmetric matrix $\Theta  \Theta ^{ *}$ is 
invertible, and
$$
\Psi := (\Theta \Theta^{ *})^{-\frac{1}{2}}
$$ 
is a bounded function of $(t,y)$.
\end{remark}

\begin{assumption}                            \label{asm3.2.18}
The random vector $z_{0}=(x_{0},y_{0})$ is 
independent of the process
$w_{t}$.  

\end{assumption}

Next, we introduce a few more notation.
Let
$$
\Psi_{t}=\Psi(t,y_{t}),\quad\Theta_{t}=\Theta(t,y_{t}),
$$
$$
a_{t}(x) =\frac{1}{2}\theta \theta^{*}(t,x,y_{t}),
\quad b_{t}(x)=b(t,x,y_{t}),\quad B_{t}(x)=B(t,x,y_{t}),
$$
$$
 \theta_{t}(x)=\theta(t,x,y_{t}),\quad
\sigma_{t}(x) =-\theta_{t}(x) \Theta^{*}_{t}\Psi_{t},\quad
\beta_{t}(x) =\Psi_{t}B_{t}(x).
$$
The columns $\sigma^{k}_{t}(x)$, $k=1,...,d_{1}$,
of the matrix-valued function $\sigma_{t}(x)$
will play the role of $\sigma^{k}_{t}(x)$, $k=1,...,d_{1}$,
 in the setting of
Section \ref{section 7.16.1}.
For $i=1 ,...,
d_{1}$ the vector-functions $\sigma^{d_{1}+i}_{t}(x)$   are defined 
as the $i$th columns of 
$$
\theta _{t}(x)-\theta _{t}(x)
\Theta^{*}_{t}
\Psi^{2}_{t}\Theta_{t} 
$$
(so that $d_{2}=d_{1}$ in the notation of Section 
\ref{section 7.16.1}, $\sigma^{0}_{t}$ will be introduced
later).
Observe that if $d'=d_{1}$,
then $\Theta_{t}$ takes values in the set of square
$d_{1}\times d_{1}$ matrices and owing to
Assumption \ref{assumption 8.21.1} is nondegenerate.
It follows that   $\Theta^{*}
\Psi^{2}\Theta$ is the identity operator 
for all $(t,y)$. In that case
$\sigma^{d_{1}+i}_{t}\equiv0$, $i=1,...,d_{1}$,
and there is no hope to get any smoothness
of the posterior distribution of $x_{t}$
unless the initial distribution has a smooth density.
Therefore, we impose the following.

\begin{assumption}
                                     \label{assumption 8.22.1}
We have $d'<d_{1}$ and, for an $S\in[0,T)$  
and $R_{0}\in(0,\infty)$,
Assumption \ref{assumption 8.2.1} is satisfied.
\end{assumption}

In a very popular so-called triangular scheme in which
$\Theta=(0,\hat{\Theta})$, where $\hat{\Theta}$ is a 
nondegenerate square $d'\times d'$-matrix valued function,
one can easily check that for $i=1 ,...,
d_{1}-d'$
$$
\sigma^{d_{1}+i}_{t}(x)=\theta^{i}_{t}(x)
$$
and $\sigma^{d_{1}+i}_{t}(x)\equiv0$ for $i=d_{1}-d'+1,...,d_{1}$.

Introduce $\cF^{y}_{t}$ as the completion of $\sigma\{y_{s}:s\leq t\}$
and denote by $P_{0}$ the regular version of the conditional
distribution of $x_{0}$ given $y_{0}$.

\begin{theorem}
                                            \label{theorem 8.22.1}
Take an $s_{0}\in(S,T]$.
Let $n$ be a negative integer such that $n<-1-d/2$.
Then there exists a function $\pi $ of class 
$H^{n+1}_{2}(T,\cF^{y}_{\cdot})$,
such that 

(i)
$\pi_{t}$ is a continuous $H^{n}_{2}$-valued
function on $[0,T]$, 

(ii) 
(a.s.) $\pi_{t}(x)$ is  infinitely differentiable 
with respect to $x$ for
$(t,x)\in [s_{0},T]\times B_{R_{0}}$ and each 
derivative is a continuous function
in $[s_{0},T]\times B_{R_{0}}$,

(iii) if $S=0$ and
the closed
support of  $ P_{0}$ is a  subset of $ B_{R_{0}}$, then (a.s.)  
$\pi_{t}(x)$  is
 infinitely differentiable 
with respect to $x$ for
$(t,x)\in ([s_{0},T]\times \bR^{d})
\cup([0,T]\times B_{R_{0}}^{c})$ and each 
derivative is a bounded continuous function
in $([s_{0},T]\times \bR^{d})\cup([0,T]\times B_{R_{0}}^{c})$, 

(iv) for any $f\in C^{\infty}_{0}(\bR^{d})$
and $t\in[0,T]$ with probability one
$$
(\pi_{t},f)=E\{f(x_{t})\mid \cF^{y}_{t}\}.
$$

\end{theorem}

Before proving the theorem we prove the following. 

\begin{lemma}
                                           \label{lemma 8.25.1}
There exists an $H^{n+1}_{2}$-valued
weakly continuous $\cF^{y}_{t}$-adapted process
$\pi_{t}$ such that assertion (iv) of Theorem \ref{theorem 8.22.1}
holds.

\end{lemma}

Proof. Take an  $f\in C^{\infty}_{0}(\bR^{d})$.
 By the famous Fujisaki-Kallianpur-Kunita theorem
the process $E\{f(x_{t})\mid \cF^{y}_{t}\}$ 
has a continuous modification, which we denote by $P_{t}(f)$.
Then a well-known procedure (see, for instance, Chapter 5, \S 3.3
\cite{Ro} or the Appendix in \cite{KR_78}) allows us to  
further modify, if necessary,  $P_{t}(f)$, so that
the new modification for which we use the same notation

(i)
is continuous in $t$ and $\cF^{y}_{t}$-adapted,

(ii) for any $\omega\in\Omega$, $t\geq 0$,
and any nonnegative $f\in C^{\infty}_{0}(\bR^{d})$ we have
$$
0\leq P_{t}(f)\leq\sup f,
$$

(iii)
  for any
$\omega\in\Omega$, $t\geq 0$,  $f,g\in C^{\infty}_{0}(\bR^{d})$,
and numbers $\alpha,\beta$ we have
$$
P_{t}(\alpha f+\beta g)=\alpha P_{t}(f)+\beta P_{t}(g).
$$

Then by Riesz-Markov theorem, there exists a measure-valued function
$P_{t}(dx)$ with $P_{t}(\bR^{d})\leq1$ such that
$$
P_{t}(f)=\int_{\bR^{d}}f(x)\,P_{t}(dx)
$$
for any $\omega\in\Omega$, $t\geq 0$,
and   $f\in C^{\infty}_{0}(\bR^{d})$.

By recalling that finite measures on $\bR^{d}$ belong to
$H^{n+1}_{2}$, we identify $P_{t}(dx)$ with a generalized 
function $\pi_{t}\in H^{n+1}_{2}$. Observe that
$$
\|\pi_{t}\|_{H^{n+1}_{2}}\leq NP_{t}(\bR^{d}),
$$
where $N$ is the embedding constant. This and the continuity
of $(\pi_{t},f)$ for $f\in C^{\infty}_{0}(\bR^{d})$ shows that
$\pi_{t}$ is 
weakly continuous as a $H^{n+1}_{2}$-valued function.
All other assertions of the lemma follow from the above.
The lemma is proved.

\begin{remark}
                                                    \label{remark 8.25.1}
If $f(t,x,y)$ is a Borel bounded function,
then for any $t\geq0$ with probability one
\begin{equation}
                                                    \label{8.25.3}
E\{f(t,x_{t},y_{t})\mid \cF^{y}_{t}\}=
\int_{\bR^{d}}f(t,x,y_{t})\,P_{t}(dx) ,
\end{equation}
where the right-hand side is a predictable function
with respect to $\{\cF^{y}_{t}\}$.

Indeed, we have seen that
$$
E\{f( x_{t} )\mid \cF^{y}_{t}\}=\int_{\bR^{d}}f(x)\,P_{t}(dx)
$$
(a.s.) for any $t\geq0$ and $f\in C^{\infty}_{0}(\bR^{d})$,
where the right-hand side is $\cF^{y}_{t}$-predictable.
This implies our claim   in an absolutely standard way.
 \end{remark}

{\bf Proof of Theorem \ref{theorem 8.22.1}}. 
Set
\begin{equation}                             \label{eq3.2.19.2}
L_{t}( x) = a^{ij}_{t}(x)D_{ij} +
 b^{i}_{t}(x)D_{i}\,,
\end{equation}
$$
L^{*}_{t}( x)u_{t}(x) = 
D_{ij} (
a^{ij}_{t}( x)u_{t}(x) )
- D_{i}(b^{i}_{t}( x)u_{t}(x) ),
$$
\begin{equation}                     \label{1.27.8}
\Lambda^{k }_{t}( x)u_{t}(x)  =\beta^{k}_{t}( x)u_{t}(x) +
 \sigma^{ik}_{t}( x) D_{i}u_{t}(x),
\end{equation}
$$
\Lambda^{k*}_{t}( x)u_{t}(x)  =\beta^{k}_{t}( x)u_{t}(x) -
D_{i}(\sigma^{ik}_{t}( x) u_{t}(x) )
$$ 
where $t\in[0,\infty)$,   $x\in \bR^{d}$, $k=1,...,d_{1} $,
and as above
 we use the summation convention over all ``reasonable''
values of repeated indices,
so that the summation in \eqref{eq3.2.19.2} and
\eqref{1.27.8}   is done for $i,j=1,...,d$ (whereas in
 \eqref{8.25.4}  for $k=1,...,d_{1}$).

By the  Fujisaki-Kallianpur-Kunita theorem, Lemma \ref{lemma 8.25.1},
and Remark \ref{remark 8.25.1} for any
$f\in C^{\infty}_{0}(\bR^{d})$ with probability one
for all $t\geq 0$
\begin{equation}
                                                    \label{8.25.4}
 (\pi_{t},f)= (\pi_{0},f)+\int_{0}^{t}
(\pi_{s},L_{s}f)\,ds+ \int_{0}^{t}(\pi_{s},\Lambda^{k}_{s}f-
\bar{\beta}^{k}
_{s}f)\,d\hat{w}^{k}_{s},
\end{equation}
where $(\hat{w}^{1}_{t},...,\hat{w}^{d_{1}}_{t})$ is a standard
Wiener process with respect to the filtration $\{\cF_{t}^{y}\}$
and
$$
\bar{\beta}_{t}=\int_{\bR^{d}}\beta(t,x,y_{t})\,P_{t}(dx).
$$

By a classical result about the It\^o formula for the squared
norm of Banach space-valued processes, equation \eqref{8.25.4}
implies that (a.s.) $\pi_{t}$
is a continuous $H^{n}_{2}$-valued process
(see Remark \ref{remark 7.31.2}). Next,  
elementary computations show that
$$
a_{t}=(1/2)\sum_{k=1}^{2d_{1}}\sigma^{k}_{t}\sigma^{k*}_{t},
\quad
L^{*}_{t} =(1/2)\sum_{k=1}^{2d_{1}}L_{\sigma^{k}_{t}}^{2} 
+L_{\sigma^{0}_{t}}+c_{t},
$$
where
$$
\sigma^{0}_{t}=(1/2)L_{\sigma^{k}_{t}}\sigma^{k}_{t}
+(\div \sigma^{k}_{t})\sigma^{k}_{t}-b_{t},\quad
c_{t}=D_{ij}a^{ij}_{t}-\div b_{t}
$$
and for $k=1,...,d_{1}$
$$
\Lambda^{k*}_{t}-\bar{\beta}^{k}=L_{\sigma^{k}_{t}}+\nu^{k}_{t},
$$
where 
$$
\nu^{k}_{t}=\beta^{k}_{t}-\bar{\beta}^{k}_{t}-\div \sigma^{k}_{t}. 
$$

After that the remaining assertion (ii) and (iii) of the present theorem
follows directly from Theorems \ref{theorem 7.22.1} 
and \ref{theorem 8.18.1}. The theorem
is proved.

In the following two theorems we, actually, speak about
$\pi_{t}(y,x)$ which is defined as $\pi_{t}(x)$ when
$P_{0}$ is the $\delta$-function concentrated at $y$.
These theorems are just direct consequences
of Theorems \ref{theorem 8.26.1} and \ref{theorem 9.11.1}
and of what was said in the proof of Theorem
\ref{theorem 8.22.1}.

\begin{theorem}
                                        \label{theorem 9.7.1}
Assume that $S=0$ and take $r\in(0,R_{0})$  
and $s_{0}\in(0,T)$. Then there exists a
nonnegative function $p_{t}(x,y)=p_{t}(\omega,y,x)$
defined for 
$$
(\omega,t,y,x)\in\Omega\times(0,T]\times B_{r}\times \bR^{d}
$$
such that

(i) it is infinitely differentiable with respect to $y$
in $B_{r}$,
 each of its $y$-derivative
of any order is a bounded function of $(t,y,x)\in[s_{0},T]
\times B_{r}\times \bR^{d}$, and the functions
$$
\int_{B_{r}}|D^{\alpha}_{y}p_{t}(y,x)|^{2}\,dy
$$
are bounded   on $[s_{0},T]
\times B_{r}$
for any   $\omega\in \Omega$  and 
multi-index $\alpha$,

(ii) for   any probability distribution
$P_{0}$, which is concentrated in $B_{r}$, 
with probability one we have
$$
\pi_{t}(x)=\int_{B_{r}} 
p_{t}( y,x)\,P_{0}(dy)
$$
for all $(t,x)\in(0,T]\times\bR^{d}$,

(iii) for any $y\in B_{r}$ with probability one
for any $t\in(0,T]$ and 
multi-index $\alpha$ the function
$D^{\alpha}_{y}p_{t}(y,x)$ is infinitely
differentiable with respect to $x$ and each of its $x$-derivative
of any order is a bounded  
function of $(t, x)\in[s_{0},T]
 \times \bR^{d}$,

(iv) the function $p_{t}(y,x)$ is infinitely differentiable
with respect to $(x,y)\in B_{r}\times B_{r}$ and each its derivative
of any order is a bounded  
 function of $(t,y,x)\in[s_{0},T]\times B_{r}^{2}$
for any  $\omega$.

\end{theorem}

We have somewhat better properties of $p_{t}(x,y)$
for a special class of filtering problems.

\begin{theorem}
                                        \label{theorem 9.7.2}
Assume that $S>0$, take $r\in(0,R_{0})$, 
$s_{0}\in (S,T)$
and assume that there exists an $R\in(0,\infty)$
such that, if $|x|\geq R$, then
 $\theta(t,x,y)=0$ and $B(t,x,y)$ is independent
of $x$.   Then there exists a
nonnegative function $p_{t}(x,y)=p_{t}(\omega,y,x)$
defined for 
$$
(\omega,t,y,x)\in\Omega\times(S,T]\times \bR^{d}\times B_{r}
$$
such that

(i)   it is infinitely differentiable
with respect to $(x,y)\in \bR^{d}\times B_{r}$, each its derivative
of any order is a bounded
 function of $(t,y,x)\in[s_{0},T]\times \bR^{d}\times B_{r}$
and the functions $\|p_{t}(\cdot,x)\|_{l}$
are bounded on $[s_{0},T]
\times B_{r}$
for any  $\omega\in \Omega$  and 
$l\geq0$,

(ii)   
with probability one we have
$$
\pi_{t}(x)=\int_{\bR^{d}} 
p_{t}( y,x)\,P_{0}(dy)
$$
for all $(t,x)\in(S,T]\times B_{r}$.

\end{theorem}

\mysection{Fundamental solutions of SPDEs}
                                         \label{section 9.10.2}

Here we continue our investigation of solutions
of general equations  \eqref{6.28.1} under the assumptions
stated in Section \ref{section 8.26.1}. We also assume that
$$
g^{k}\equiv0,\quad
  k=1,...,d_{1} .
$$
Here again $u$ is the normal solution from Theorem \ref{theorem 7.24.1}
and $n$ is from Assumption \ref{assumption 6.28.1}.

\begin{lemma}
                                            \label{lemma 8.2.1}
Take an $R\in[0,\infty)$, $s_{0}\in(S,T)$, $r\in(0,R_{0})$,
and assume that $\sigma^{k}_{t} =0$ and $\nu^{k}_{t} =0$ outside $B_{R}$ for
any $k=1,...,d_{1}$, $t$, and $\omega$. 
 Then,
for any  $l$ and multi-index $\alpha$, such that
\begin{equation}
                                                 \label{7.16.01} 
2(l-|\alpha|-2)>d+1,
\end{equation}
with probability one
there exists a (random, finite)
constant $N$ independent of $u$ such that
\begin{equation}
                                                 \label{7.22.02}
\sup_{(t,x)\in[s_{0},T]\times B_{r}}
|D^{\alpha}u_{t}(x)|^{2}\leq N\|u_{0}\|_{n}^{2}+
N\int_{0}^{T}\|f_{t} \|_{l}^{2}\,dt.
\end{equation}
\end{lemma}

 Proof. Obviously,
it suffices to concentrate on $n\leq l$.
Take a $\zeta\in C^{\infty}_{0}(B_{R_{0}})$ such that $\zeta=1$
in a neighborhood of $\bar{B}_{r}$. Then by Theorem \ref{theorem 7.22.1}  
for any  $l$ and multi-index $\alpha$ satisfying \eqref{7.16.01}
and $n\in\bR$,
with probability one,
there exists a (random, finite)
constant $N$ independent of $u$ such that  
\begin{equation}
                                                 \label{7.22.002}
\sup_{(t,x)\in[S,T]\times B_{r}}
|D^{\alpha}u_{t}(x)|^{2}\leq N
\int_{S}^{T}\big[\|f_{t}\zeta \|_{l}^{2}
+\| u_{t}  
\zeta \|_{n}^{2}\big]\,dt.
\end{equation}
Owing to Theorem \ref{theorem 8.1.1}
$$
\| u_{t}  
\zeta \|_{n}  ^{2}\leq N\| u_{t}  
  \|_{n }^{2}\leq N\|u_{0}\|_{n}^{2}+
N\int_{0}^{T}\|f_{s} \|_{n}^{2}\,ds
$$
and this proves the lemma.

\begin{lemma}
                                            \label{lemma 8.14.2}
In the setting of of Lemma \ref{lemma 8.2.1}
  suppose that $S=0$ and $u_{0} $ is infinitely differentiable
in $B_{R_{0}}$ for every $\omega$. Then
for any $\zeta\in C^{\infty}_{0}(B_{R_{0}})$ such that $\zeta=1$
in a neighborhood of $\bar{B}_{r}$, any
$l$ and multi-index $\alpha$ satisfying \eqref{7.16.01},
  with probability one,
there exists a (random, finite)
constant $N$ independent of $u$ such that
$$
\sup_{(t,x)\in[0,T]\times B_{r}}
|D^{\alpha}u_{t}(x)|^{2}\leq N\|\zeta u_{0}\|^{2}_{l+1}+
N\int_{0}^{T}\|f_{t} \|_{l}^{2}\,dt+
N\|  u_{0}\|^{2}_{n}.
$$

\end{lemma}

This lemma is derived from Theorem \ref{theorem 7.22.1}   and
Theorem \ref{theorem 8.1.1} in the same way as Lemma \ref{lemma 8.2.1}.

\begin{theorem}
                                                \label{theorem 8.14.3}
Assume that $S=0$ and take $s_{0}\in
(0,T)$ and   $0<r<r_{1}<r_{2} <R_{0}$. 
Suppose that  $u_{0}\in\cHO^{n}_{2}(B_{r})$.
Finally, let $f_{t}(x)=0$ for $|x|\geq r$.

Then
for  any
$l$ and multi-index $\alpha$ satisfying \eqref{7.16.01},
and $n\in\bR$, with probability one,
there exists a (random, finite)
constant $N$ independent of $u$ such that
\begin{equation}
                                                 \label{8.14.4}
\sup_{(t,x)\in\Gamma}
|D^{\alpha}u_{t}(x)|^{2}\leq  NJ,
\end{equation}
where 
$$
\Gamma=\big([0,T]\times B_{r_{2}}\big)\setminus
\big([0,s_{0}]\times B_{r_{1}}\big),
\quad J=\int_{0}^{T} \|f_{t}\|^{2}_{l} \,dt +N \|u_{0}\|^{2}_{n}.
$$

\end{theorem}

Proof. We take $ r_{2}<r_{3}<R_{0}$, $\zeta\in 
 C^{\infty}_{0}(B_{r_{3}})$, and $\eta\in 
 C^{\infty}_{0}(B_{r_{2}})$ such that $\zeta=1$ in a neighborhood of
 $\bar B_{r_{2}}$
and $\eta=1$  in a neighborhood of $\bar B_{r_{1}}$.

As in the proof of Theorem \ref{theorem 8.18.1}
by covering $B_{r_{2}}\setminus B_{r_{1}}$ with
appropriate balls and applying Theorem \ref{theorem 7.22.1}  
we see that
$$
\sup_{(t,x)\in[0,T]\times(B_{r_{2}}\setminus B_{r_{1}})}
|D^{\alpha}u_{t}(x)|^{2}\leq  
N\int_{0}^{T}[\|f_{t} \|_{l}^{2}+\|\zeta u_{t}\|^{2}_{n}]\,dt,
$$
where the (random) constant $N$ is independent of $u,f$.
All such constants will be denoted by $N$.
Furthermore, by the same theorem
$$
\sup_{(t,x)\in[s_{0},T]\times B_{r_{2}} }
|D^{\alpha}u_{t}(x)|^{2}\leq  
N\int_{0}^{T}[\|f_{t} \|_{l}^{2}+\|\zeta u_{t}\|^{2}_{n}]\,dt.
$$

Obviously, we may assume that $n<l$ and then we note
 that to prove \eqref{8.14.4}, it suffices to show
that with probability one,
there exists a (random, finite)
constant $N$ independent of $u$ such that
\begin{equation}
                                                 \label{8.21.1}
\int_{0}^{T} \|\zeta u_{t}\|^{2}_{n} \,dt\leq  NJ.
\end{equation}

To prove \eqref{8.21.1} we need an auxiliary function.
Take a $\xi\in C^{\infty}_{0}(B_{2R_{0}})$
such that $\xi=1$ on $B_{R_{0}}$ and define $v_{t}$
as a normal $H^{n}_{2}(T)$-solution of the equation
$$
dv_{t}=(L_{t}v_{t}+c_{t}v_{t}+f_{t})\,dt
+(L_{\xi\sigma^{ k}_{t}}v_{t}+\xi\nu^{k}_{t}v_{t} )\,dw^{k}_{t}
$$
with initial condition $u_{0}$. 
By Theorem \ref{theorem 8.1.1}
\begin{equation}
                                        \label{8.20.1}
\sup_{t\leq T}\|v_{t}\|^{2}_{n}\leq
N J.
\end{equation}

By Theorem \ref{theorem 7.22.1}  
$$
\sup_{[0,T]}\sup_{|x|=r_{1}}|v_{s}(x)|^{2}\leq
N\int_{0}^{T}[\|f_{t} \|_{l}^{2}+\|v_{t}\|^{2}_{n}]\,dt,
$$
which after being combined with \eqref{8.20.1} shows that
$$
\sup_{[0,T]}\sup_{|x|=r_{1}}|v_{s}(x)|^{2}\leq
NJ 
$$
and
$$
\sup_{[0,t]}\sup_{|x|=r_{1}}|u_{s}(x)|^{2}
\leq 2\sup_{[0,t]}\sup_{|x|=r_{1}}|u_{s}(x)-v_{s}(x)|^{2}+NJ.
$$

By applying 
Theorem \ref{theorem 7.22.1}   to $u_{t}-v_{t}$ we conclude that
for $t\in[0,T]$
\begin{equation}
                                                      \label{8.20.3}
\sup_{[0,t]}\sup_{B_{r_{1}}}|u_{s}-v_{s}|^{2}\leq
N\int_{0}^{t}\|\eta(u_{s}-v_{s})\|^{2}_{n}\,ds.
\end{equation}
Hence,
$$
\sup_{[0,t]}\sup_{|x|=r_{1}}|u_{s}(x)|^{2}\leq 
N\int_{0}^{t}\|\eta(u_{s}-v_{s})\|^{2}_{n}\,ds
+NJ
\leq 
N\int_{0}^{t}\|\eta u_{s} \|^{2}_{n}\,ds
+NJ.
$$

Next, set
$$
P_{t}=\|u_{t}\zeta\|^{2}_{n},\quad Q_{t}=\|u_{t}\eta\|^{2}_{n}
$$
and observe that
$$
P_{t}\leq 2Q_{t}+2\|u_{t}(\zeta-\eta)\|_{n}^{2},
$$
where,
as in the proof of Theorem \ref{theorem 8.18.1}
$u_{t}(\zeta-\eta)$ is infinitely differentiable 
with respect to $x$ for
$(t,x)\in [0,T]\times \bR^{d}$ and each 
derivative is a bounded continuous function
in $[0,T]\times \bR^{d}$. In particular,
$$
\|u_{t}(\zeta-\eta)\|_{n}\leq N\sup_{[0,t]
\times \bR^{d}}|u_{s}(\zeta-\eta)(x)|
$$
$$
\leq N\sup_{[0,t]}\sup_{B^{c}_{r_{1}}}|u_{s}(x)|\leq
N\sup_{[0,t]}\sup_{|x|=r_{1}}|u_{s}(x)|,
$$
where the last inequality is taken from Theorem \ref{theorem 8.13.1}.
Hence,
\begin{equation}
                                                      \label{8.21.3}
P_{t}\leq 2Q_{t}+ 
N\int_{0}^{t}\|\eta u_{s} \|^{2}_{n}\,ds
+NJ
\leq 2Q_{t}+ 
N\int_{0}^{t}P_{s}\,ds
+NJ.
\end{equation}

Then 
$$
Q_{t}\leq 2\|(u_{t}-v_{t})\eta\|^{2}_{n}+2\|v_{t}\|^{2}_{n}
\leq N\sup_{ B_{r_{2}}}| u_{t}-v_{t} |^{2}+NJ
$$
and similarly to \eqref{8.20.3}
$$
\sup_{ B_{r_{2}}}| u_{t}-v_{t} |^{2}\leq 
N\int_{0}^{t}\|\zeta(u_{s}-v_{s})\|^{2}_{n}\,ds
\leq 
N\int_{0}^{t}P_{s}\,ds+NJ.
$$
By coming back to \eqref{8.21.3} and using Gronwall's inequality
we conclude $P_{t}\leq NJ$ and this proves \eqref{8.21.1}
and the theorem.

In th following theorem we prove the existence
of a kernel for our SPDE. It is worth
drawing the reader's attention to the fact that
no continuity with respect to 
 $(t,y,x)$ is claimed in (i) and (iv)
and no continuity with respect to $y$ is
claimed in (iii).
 
\begin{theorem} 
                                    \label{theorem 8.26.1}
Assume that $S=0$ and take $r\in(0,R_{0})$
and $s_{0}\in(0,T)$. Then there exists a
nonnegative function $p_{t}(x,y)=p_{t}(\omega,y,x)$
defined for 
$$
(\omega,t,y,x)\in\Omega\times(0,T]\times B_{r}\times \bR^{d}
$$
such that  

(i) it is infinitely differentiable with respect to $y$
in $B_{r}$,
 each of its $y$-derivative
of any order is a bounded
  function of $(t,y,x)\in[s_{0},T]
\times B_{r}\times \bR^{d}$, and the functions
\begin{equation}
                                               \label{9.12.1}
\int_{B_{r}}|D^{\alpha}_{y}p_{t}(y,x)|^{2}\,dy
\end{equation}
are bounded on $[s_{0},T]
\times B_{r}$
for any   $\omega\in \Omega$  and 
multi-index $\alpha$,

(ii) if $u_{0}\in \cHO^{n}_{2}(B_{r})$
and $f\equiv g^{k}\equiv0$,
$k=1,...,d_{1}$, then
with probability one $u_{t}(x)$ coincides (as a generalized function
with respect to $x$) with
\begin{equation}
                                                    \label{8.28.2}
(u_{0},p_{t}(\cdot,x))=\int_{B_{r}}u_{0}(y)
p_{t}( y,x)\,dy
\end{equation}
for all $(t,x)\in(0,T]\times\bR^{d}$.

(iii) for any $y\in B_{r}$ with probability one
for any $t\in(0,T]$ and 
multi-index $\alpha$ the function
$D^{\alpha}_{y}p_{t}(y,x)$ is infinitely
differentiable with respect to $x$ and each of its $x$-derivative
of any order is a bounded
and continuous function of $(t, x)\in[s_{0},T]
 \times \bR^{d}$,

(iv) the function $p_{t}(y,x)$ is infinitely differentiable
with respect to $(x,y)\in B_{r}\times B_{r}$ and each its derivative
of any order is a bounded  
 function of $(t,y,x)\in[s_{0},T]\times B_{r}^{2}$
for any  $\omega\in \Omega$.
\end{theorem}

Proof. Let $\cZ_{+}=\{\zeta_{1},\zeta_{2},
...\}$ be a countable subset
of $C^{\infty}_{0}(B_{r})$ consisting of nonnegative
functions such that 
 $\cZ_{+}$ is dense in the sup-norm in the subset of $C^{\infty}_{0}(B_{r})$ 
consisting of nonnegative functions.
Also set $\cZ=\cZ_{+}-\cZ_{+}$.

For $\zeta\in \cHO^{-l}_{2}(B_{r})$
 denote by $u_{t}[\zeta](x)$ the normal solution
 of \eqref{6.28.1} with initial condition
$\zeta$ and with $f\equiv g^{k}\equiv0$,
$k=1,...,d_{1}$. By Theorems \ref{theorem 8.14.3} and
\ref{theorem 8.13.1} for any  $l\geq0$ and  $t\in(0,T]$
with probability one there exists a constant $N_{t}(l)$
such that for any $x\in\bR^{d}$ and $\zeta\in \cHO^{-l}_{2}(B_{r})$ 
we have
\begin{equation}
                                                \label{8.27.1}
|u_{t}[\zeta](x)|\leq N_{t}(l)
\|\zeta\|_{-l},
\end{equation}
and if $\zeta\in\cZ_{+}$ then by the maximum principle
\begin{equation}
                                                \label{8.27.2}
 u_{t}[\zeta](x)\geq0.
\end{equation}

Furthermore, owing to uniqueness for any 
collection of rational numbers $p_{1},p_{2},...$
only finite number of which is different
from zero with probability one
\begin{equation}
                                                \label{8.27.3}
\big|\sum_{i}p_{i}u_{t}[\zeta_{i}](x)\big|=
\big|u_{t}[\sum_{i}p_{i}\zeta_{i}](x)\big|
\leq N_{t}(l)\big\|\sum_{i}p_{i}\zeta_{i}\big\|_{-l}.
\end{equation}

Obviously, one can choose $N_{t}(l)$ so that it is
a monotone function of $t$. Owing to this and the fact that
other expressions entering \eqref{8.27.3}
and \eqref{8.27.2} are continuous in $(t,x)$
(a.s.) and the fact that the set of rational numbers is
countable,  there is a set $\Omega'$ of full probability
such that \eqref{8.27.3} is satisfied
   for any $l$, any
collection of rational numbers $p_{1},p_{2},...$,
only finite number of which is different
from zero,   any $\omega\in\Omega'$, and any
$(t,x)\in(0,T]\times\bR^{d}$, and  \eqref{8.27.2}
is satisfied  for any $\omega\in\Omega'$ and any
$(t,x)\in(0,T]\times\bR^{d}$ and $\zeta\in\cZ$.
By setting
$u_{t}[\zeta](x)=0$ outside $\Omega'$, we may assume that
$\Omega'=\Omega$.

By a theorem of Hahn (see, for instance, Section II.5  
of \cite{DS}), for any $x\in\bR^{d}$, $t\in(0,T]$, $l\geq0$,
and $\omega$ there exists a linear functional $Q_{t}[\cdot](x)$
on $H^{-l}_{2}$ such that
$$
Q_{t}[\zeta_{i}](x)=u_{t}[\zeta_{i}](x),\quad\forall i,
\quad\|Q_{t}[\cdot](x)\|\leq N_{t}(l).
$$
The general form of linear functionals $Q$ on $H^{-l}_{2}$
is well known and for smooth elements   $f\in H^{-l}_{2}$
it is given by
$$
Q(f)=\int_{\bR^{d}}f(x)p(x)\,dx,
$$
where $p\in H^{l}_{2}$ and $\|p\|_{l}=\|Q\| $.

Hence, for any $x\in\bR^{d}$, $t\in(0,T]$, $l\geq0$,
and $\omega$ there exists 
a  function   $p_{t}(\cdot,x)
\in H^{l}_{2}$  such that
\begin{equation}
                                                \label{8.27.4}
\|p_{t}(\cdot,x)\|_{l}\leq N_{t}(l),
\quad \int_{\bR^{d}}\zeta_{i}(y)p_{t}(y,x)\,dy=u_{t}[\zeta_{i}](x)
\quad\forall i.
\end{equation}

In principle $p_{t}(y,x)$ is not unique, it might be
  changed for $y\not\in B_{r}$. In addition it may depend
on $n$.
To choose
a better representative, without loosing generality 
 we take $l=2m$, where $m=1,2,...$,
so that
$$
\|u\|_{l}=\|(1-\Delta)^{m}u\|_{0}
$$
and then set $p_{t}(y,x)=0$ for $y\not\in B_{r}$. Then the new
$p_{t}(\cdot,x)$ for which we keep the same notation will
satisfy the second relation  in \eqref{8.27.4}, hence
will be independent of $l$. Also for thus modified
$p_{t}(y,x)$ the functions \eqref{9.12.1}
 are bounded on $[s_{0},T]
\times B_{r}$
for any $s_{0}\in(0,T]$, $\omega\in \Omega$, and 
multi-index $\alpha$ since this was true for the initial
$p_{t}(y,x)$ even with   $\bR^{d}$ in place of $B_{r}$ 
in \eqref{9.12.1}.

 Since this holds for any multi-index,  applying 
embedding theorems we obtain assertion (i). 
The second relation   in \eqref{8.27.4} along with 
\eqref{8.27.3} and the choice
of $\{\zeta_{i}\}$ implies that $p_{t}(y,x)\geq0$
for any $x\in\bR^{d}$, $y\in B_{r}$, $t\in(0,T]$,  
and $\omega$.

Estimate \eqref{8.27.1}, equations in \eqref{8.27.4},
and the fact that $\cZ $
is dense in $\cHO^{ n}_{2}(B_{r})$ implies 
assertion (ii).
The expressions in \eqref{8.28.2} are (s.s.)
 infinitely differentiable
with respect to $x$ and each of the derivatives is
bounded and continuous in $[s_{0},T]\times\bR^{d}$
by Theorem \ref{theorem 8.18.1}.  
This proves assertion (iii)
if we set $u_{0}=D^{\alpha}\delta_{y}$,
where $\delta_{y}$ is the delta-function concentrated
at $y\in B_{r}$.

To prove assertion (iv) observe that, as above,
owing to Theorem \ref{theorem 8.14.3}, we can modify,
if necessary, the functions
$u_{t}[\zeta ](x)$ for $\zeta\in\cZ $ in such a way that
for any $\omega$ and $t\in(0,T]$ they will be infinitely differentiable
 with respect to $x$ in $B_{r}$ and for any
multi-index $\alpha$ and $l>0$ satisfy
\begin{equation}
                                                     \label{8.29.1}
\sup_{B_{r}}|D^{\alpha}u_{t}[\zeta](x)|\leq
N_{t}(|\alpha|,l)\|\zeta\|_{-l},
\end{equation}
where $N_{t}(|\alpha|,l)$ are independent of $\zeta$.
In particular, for $x,x',x''\in B_{r}$ and $\zeta,\zeta',\zeta''
\in \cZ $
$$
|u_{t}[\zeta](x')-u_{t}[\zeta](x')|\leq
N_{t}(1,l)\|\zeta\|_{-l}|x'-x''|,
$$
\begin{equation}
                                                     \label{8.29.2}
|u_{t}[\zeta'](x)-u_{t}[\zeta''](x)|\leq N_{t}(0,l)
\|\zeta'-\zeta''\|_{-l}.
\end{equation}

Since $\cZ $ is dense in $\cHO^{-l}_{2}(B_{r})$ and 
\eqref{8.27.4} holds, estimates  \eqref{8.29.2}
hold for all $\zeta,\zeta',\zeta''
\in \cHO^{-l}_{2}(B_{r})$, which after taking an appropriate
$l$ and $\zeta,\zeta',\zeta''$ as $\delta$-functions show
that $p_{t}(x,y)$ is a continuous function of $(x,y)$
on $B_{r}^{2}$.

Estimate \eqref{8.29.1} also implies that
the generalized derivatives
of $(\zeta,p_{t}(\cdot,x))$ of order $|\alpha|$
are bounded on $B_{r}$ by the right-hand side of \eqref{8.29.1}
for any $\zeta\in\cHO^{-l}_{2}(B_{r})$.
Since this holds for any $\alpha$
and $(\zeta,p_{t}(\cdot,x))$ is continuous with respect to $x\in B_{r}$
(owing to \eqref{8.29.2}),
 the generalized derivatives are, actually,
usual ones, which are bounded and continuos.

By taking $\zeta=D^{\beta}_{y}\delta_{y}$ in \eqref{8.29.1}, 
we find that  the usual (and hence generalized)
 functions $D^{\alpha}_{x}[
D^{\beta}_{y}p_{t}(y,x)]$ are bounded for $(x,y)\in B_{r}^{2}$.
It follows that $p_{t}(y,x)$ admits a modification
with respect to $(y,x)$ which is infinitely differentiable.
However, the modification coincides with $p_{t}(y,x)$
on $B_{r}^{2}$ since $p_{t}(y,x)$ is continuous there.

The asserted boundedness of the derivatives of $p_{t}(y,x)$
is easily derived from the above argument.
The theorem is proved.

In the following theorem we assert the regularity
of 
$p_{t}(y,x)$ not only for   $y\in B_{r}$
but for all $y\in\bR^{d}$ albeit for $x\in B_{r}$,
the latter being of course inevitable.
 But the result is proved
under a somewhat restrictive assumption. This assumption
arose because of out inability to control
$\|u_{t}\|_{n}$ through $\|u_{0}\|_{n}$ times
 a (random) constant independent of $u_{0}$.

\begin{theorem}
                                        \label{theorem 9.11.1}
Suppose that $S>0$,
take an $R\in[0,\infty)$, $r\in(0,R_{0})$,  
$s_{0}\in(S,T)$,
and assume that $\sigma^{k}_{t} =0$ and $\nu^{k}_{t} =0$ outside $B_{R}$ for
any $k=1,...,d_{1}$, $t$, and $\omega$.  Then there exists a
nonnegative function $p_{t}(x,y)=p_{t}(\omega,y,x)$
defined for 
$$
(\omega,t,y,x)\in\Omega\times(S,T]\times \bR^{d}
\times B  _{r}
$$
such that

(i)  it is infinitely differentiable
with respect to $(x,y)\in \bR^{d}\times B_{r}$, each its derivative
of any order is a bounded
 function of $(t,y,x)\in[s_{0},T]\times \bR^{d}\times B_{r}$
and the functions $\|p_{t}(\cdot,x)\|_{l}$
are bounded on $[s_{0},T]
\times B_{r}$
for any  $\omega\in \Omega$  and 
$l\geq0$,

(ii) if $f\equiv g^{k}\equiv0$,
$k=1,...,d_{1}$, then
\begin{equation}
                                                  \label{9.11.2}
u_{t}(x)=\int_{\bR^{d}} 
p_{t}( y,x)u_{0}(y) \,dy 
\end{equation}
for all $(t,x)\in(S,T]\times B_{r}$.

\end{theorem}

Proof. It suffices to repeat the proof of Theorem 
\ref{theorem 8.26.1} using the same notation as there,
replacing $\cZ$ with a countable subset
of $C^{\infty}_{0}(\bR^{d})$, and using
  Lemma \ref{lemma 8.2.1} to conclude that  
 for any 
and $l\geq0$ and $t\in(S,T]$
with probability one there exists a constant $N_{t}(l)$ 
such that, for any   $\zeta\in \cHO^{-l}_{2} $  and
 $x\in B_{r}$, estimate \eqref{8.27.1}
holds.
After that the proof goes almost exactly the same way
as that of Theorem 
\ref{theorem 8.26.1}. The theorem is proved.

\end{document}